# Sampled-Data Observers for Delay Systems and Hyperbolic PDE-ODE Loops

**Tarek Ahmed-Ali[*], Iasson Karafyllis[**] and Fouad Giri[*]**

[*]Normandie UNIV, UNICAEN, ENSICAEN, LAC, 14000 Caen, France
email: tarek.ahmed-ali@ensicaen.fr; fouad.giri@unicaen.fr

[**]Dept. of Mathematics, National Technical University of Athens,
Zografou Campus, 15780, Athens, Greece,
email: iasonkar@central.ntua.gr; iasonkaraf@gmail.com

**Abstract**
This paper studies the problem of designing sampled-data observers and observer-based, sampled-data, output feedback stabilizers for systems with both discrete and distributed, state and output time-delays. The obtained results can be applied to time delay systems of strict-feedback structure, transport Partial Differential Equations (PDEs) with nonlocal terms, and feedback interconnections of Ordinary Differential Equations with a transport PDE. The proposed design approach consists in exploiting an existing observer, which features robust exponential convergence of the error when continuous-time output measurements are available. The observer is then modified, mainly by adding an inter-sample output predictor, to compensate for the effect of data-sampling. Using Lyapunov stability tools and small-gain analysis, we show that robust exponential stability of the error is preserved, provided the sampling period is not too large. The general result is illustrated with different examples including state observation and output-feedback stabilization of a chemical reactor.

## 1. Introduction

With the growing penetration of network technology in control systems, the compensation of the effects of time-delay has become a major issue in control theory [19,25,26,27,33]. A great deal of interest has recently been paid to the problem of designing state observers for linear and nonlinear systems with measurement delays. The dominant design approach consists in starting with the design of an exponentially convergent observer for the delay-free system, which is described by Ordinary Differential Equations (ODEs), and modifying it mainly by adding predictors: static predictors (see [27]) or dynamic (chain) predictors (see [1,4,12,14,15,30]). In parallel to this research activity, which takes into account the time-delay explicitly in the model, a separate activity, based on Partial Differential Equations (PDEs) has been initiated in [31]. This consists in modeling time-delays by means of first-order hyperbolic PDEs, leading to a representation of the delayed system in the form of an ODE-PDE cascade (see also the recent work [3], where a PDE-based chain-observer is constructed).

Most existing results on observer design for delayed systems have been established assuming the measurement delay to be of discrete nature. So far, only a few studies have investigated the case of distributed measurement time-delays. The PDE-based observer developed in [10] and the recent observer developed in [7], are notable exceptions.

The nowadays-digital implementation of observers entails sampling in time of all system signals needed by the observer. Consequently, the design of sampled-data observers is a major issue. Sampled-data observers for ODE systems can be classified in four main categories:

1) observers where data-sampling is accounted for through a standard Zero-Order-Hold (ZOH) sampling of the output estimation error; see for example [2,37],
2) observers designed on approximate discrete-time models (see [8,13]),
3) continuous-discrete time observers where correction terms are employed at the sampling times; see for instance [34], and
4) sampled-data observers, where the time-varying delay effect (caused by output sampling) is compensated by using inter-sample output predictors; see [23].

The use of inter-sample output predictors was extended to systems with asynchronous measurements (see [32]) and systems described by parabolic PDEs (see [29]).

The combination of time-delay and data-sampling effects necessarily makes the problem of observer design more complex. Indeed, not only data-sampling introduce a time-varying delay but it also entails information lost. The case of discrete measurement delays, in conjunction with data sampling, has been investigated in [4,5,25,27,37]. Results on observer-based output feedback stabilization of delay systems with sampled measurements have been recently given in [16,36] (but see also the case of state feedback in [35]).

In the present work, we extend for the first time the use of inter-sample predictors to the case of time-delay systems with state and output (discrete and/or distributed) measurement delays. Moreover, we provide observer-based output feedback stabilization results for delay systems with sampled measurements under appropriate assumptions. More specifically, we consider nonlinear time-delay systems of the form:

$$\dot{x} = f(x_t, u, d)$$
$$y = h(x_t) + \xi \quad (1.1)$$
$$(x, u, d) \in \Re^n \times U \times D, \, y, \xi \in \Re^k$$

where $U \subseteq \Re^m$, $D \subseteq \Re^q$ are convex sets with $0 \in U$, $0 \in D$, $f: C^0([-r, 0]; \Re^n) \times U \times D \to \Re^n$, $h: C^0([-r, 0]; \Re^n) \to \Re^k$ are continuous mappings with $f(0,0,0) = 0$, $h(0) = 0$. The input $u$ is assumed to be available, but the inputs $d, \xi$ are unknown and represent possible modeling errors and measurement noise, respectively. The proposed sampled-data observer design approach consists in starting with an existing observer, which features *robust exponential convergence* when continuous-time output measurements are available (see Definition 2.1 for the precise meaning of the phrase "robust exponential convergence"). The available observer, based on continuous-time measurements, is then modified, mainly by adding an inter-sample output predictor, to compensate for the effect of data-sampling. Using Lyapunov stability tools and small gain analysis, we show that the robust exponential stability feature is preserved, provided that the sampling period is sufficiently small (Theorem 2.2). The sampled-data observer can be used in a straightforward way for the design of observer-based output feedback stabilizers (Corollary 2.4) under certain assumptions.

The second contribution of the paper is that it provides a framework where sampled-data observer and feedback design for ODE-PDE loops can be converted to a similar problem for a time-delay system. More specifically, we consider feedback interconnections of ODEs with a first-order, hyperbolic PDE with non-local terms of the form

$$\frac{d\bar{x}}{dt}(t) = \tilde{F}(\bar{x}(t), v[t], u(t)), \text{ for } t \geq 0 \tag{1.2}$$

$$\frac{\partial v}{\partial t}(t,z) + c\frac{\partial v}{\partial z}(t,z) = a(z)v(t,z) + g(z,\bar{x}(t)) + \sum_{i=1}^{N_1} b_i(z)v(t,z_i) + \sum_{i=1}^{N_2} \beta_i(z)\int_0^1 \gamma_i(s)v(t,s)ds,$$

$$\text{for all } (t,z) \in \mathfrak{R}_+ \times [0,1] \tag{1.3}$$

$$v(t,0) = 0, \text{ for } t \geq 0 \tag{1.4}$$

where $c > 0$ is a constant (the transport speed), $z_i \in (0,1]$ ($i=1,...,N_1$), $\bar{x}(t) \in \mathfrak{R}^{\bar{n}}$, $v(t,z) \in \mathfrak{R}$ are the states, $u \in C^0(\mathfrak{R}_+; \mathfrak{R}^m)$ is an external input and the mappings $a:[0,1] \to \mathfrak{R}$, $b_i:[0,1] \to \mathfrak{R}$ ($i=1,...,N_1$), $\beta_i, \gamma_i:[0,1] \to \mathfrak{R}$ ($i=1,...,N_2$), $\tilde{F}:\mathfrak{R}^{\bar{n}} \times C^0([0,1]) \times \mathfrak{R}^m \to \mathfrak{R}^{\bar{n}}$, $g:[0,1] \times \mathfrak{R}^{\bar{n}} \to \mathfrak{R}$ are given continuous mappings. For such systems the output $y(t) = (y_1(t),...,y_k(t))^T \in \mathfrak{R}^k$ is given by equations of the following form:

$$y_j(t) = q_j^T \bar{x}(t) + \sum_{i=1}^{N_1} \bar{b}_{j,i} v(t,z_i) + \sum_{i=1}^{N_2} \bar{\beta}_{j,i} \int_0^1 \gamma_i(s)v(t,s)ds, \quad j=1,...,k \tag{1.5}$$

where $q_j \in \mathfrak{R}^{\bar{n}}$ ($j=1,...,k$) are constant vectors, $\bar{b}_{j,i} \in \mathfrak{R}$ ($j=1,...,k$, $i=1,...,N_1$) and $\bar{\beta}_{j,i} \in \mathfrak{R}$ ($j=1,...,k$, $i=1,...,N_2$) are constants. Hyperbolic PDE-ODE loops have been studied recently in [3,9,17,18,28,38]. For the class of systems (1.2), (1.3), (1.4), (1.5), we provide sampled-data observers and observer-based, sampled-data, output feedback stabilizers (Theorem 3.2 and Corollary 3.4). The usefulness of the obtained results is illustrated through the design of a sampled-data observer for a chemical reactor model proposed in [28] and the output feedback stabilization of the reactor (Theorem 3.3 and Example 3.5) with no restriction on the transport speed (or equivalently with no restriction on the induced delay). The obtained results can also be applied to the case where no ODEs are present (i.e., the case of a single first-order hyperbolic PDE with non-local terms). The case of observer design for a single first-order hyperbolic PDE was recently studied in [11]. It should also be noted that for many cases where there is an interconnection term in the boundary condition of the PDE, there is a simple transformation that can express the system to the form (1.2), (1.3), (1.4), (1.5) and consequently, the form (1.2), (1.3), (1.4), (1.5) is not restrictive.

Sampled-data observers are also provided for the class of uncertain, triangular, globally Lipschitz delay systems of the form

$$\dot{x}_i(t) = f_i(x_{1,t},...,x_{i,t},u(t)) + x_{i+1}(t) + d_i(t), \quad i=1,...,n-1$$
$$\dot{x}_n(t) = f_n(x_{1,t},...,x_{n,t},u(t)) + d_n(t) \tag{1.6}$$
$$y(t) = x_1(t)$$

where $x(t) = (x_1(t),...,x_n(t)) \in \mathfrak{R}^n$ is the state, $u(t) \in \mathfrak{R}^m$ is a known input, $d(t) = (d_1(t),...,d_n(t)) \in \mathfrak{R}^n$ is the vector of the disturbances (or unknown inputs) and $f_i : C^0([-r,0]; \mathfrak{R}^i) \times \mathfrak{R}^m \to \mathfrak{R}$ ($i=1,...,n$) with $f_i(0) = 0$ ($i=1,...,n$) are globally Lipschitz functional with $r > 0$ being the maximum delay. Again, by using the inter-sample predictor design, we are in a position to design sampled-data observers for system (1.6), no matter how large the maximum delay $r>0$ is (Theorem 4.1). The observer design is based on the high-gain observer design for ODEs, proposed in [20].

It should be noted that in all cases the results are global. Moreover, we are in a position to consider uncertain sampling schedules (i.e., the sampling times are not a priori known) and guarantee robustness with respect to measurement noise. Finally, in the absence of measurement noise and unknown disturbances, exponential convergence of the observer error is achieved. The obtained results in all cases are proved by means of a combined use of Lyapunov-Krasovskii functionals and a small-gain analysis. The fact that the proposed sampled-data observer design approach with an inter-sample output predictor can indifferently be applied to time-delay systems, to transport PDEs, and to interconnections of ODEs with a transport PDE, provides the approach with a strong unifying feature.

The structure of the paper is as follows. In Section 2 we give the general results for the time-delay systems of the form (1.1). Section 3 is devoted to the study of ODE-PDE loops of the form (1.2), (1.3), (1.4), (1.5). The results for the triangular, globally Lipschitz, time-delay case (1.6) are presented in Section 4. All the proofs are provided in Section 5. Finally, the concluding remarks of the present work are given in Section 6.

**Notation.** Throughout this paper we adopt the following notation:

* By $\Re_+$ we denote the set of non-negative real numbers. Let $u:\Re_+ \times [0,1] \to \Re$ be given. We use the notation $u[t]$ to denote the profile at certain $t \geq 0$, i.e., $(u[t])(z) = u(t,z)$ for all $z \in [0,1]$. For $u \in C^0([0,1])$ the sup norm is defined by $\|u\|_\infty = \sup_{0 \leq z \leq 1}(|u(z)|)$.

* Let $S \subseteq \Re^n$ be an open set and let $A \subseteq \Re^n$ be a set that satisfies $S \subseteq A \subseteq cl(S)$. By $C^0(A;\Omega)$, we denote the class of continuous functions on $A$, which take values in $\Omega \subseteq \Re^m$. By $C^k(A;\Omega)$, where $k \geq 1$ is an integer, we denote the class of functions on $A \subseteq \Re^n$, which takes values in $\Omega \subseteq \Re^m$ and has continuous derivatives of order $k$. In other words, the functions of class $C^k(A;\Omega)$ are the functions which have continuous derivatives of order $k$ in $S = int(A)$ that can be continued continuously to all points in $\partial S \cap A$. When $\Omega = \Re$ then we write $C^0(A)$ or $C^k(A)$.

* For a vector $x \in \Re^n$ we denote by $|x|$ its usual Euclidean norm and by $x^T$ its transpose. By $|A| := \sup\{|Ax|; x \in \Re^n, |x|=1\}$ we denote the induced norm of a matrix $A \in \Re^{m \times n}$ and $I$ denotes the identity matrix. By $B = diag(b_1,...,b_n)$ we denote the diagonal matrix $B \in \Re^{n \times n}$ with $b_1,...,b_n$ in its diagonal. For $x \in C^0([-r,0];\Re^n)$ we define $\|x\| := \max_{\theta \in [-r,0]}(|x(\theta)|)$.

* Let $x:[a-r,b) \to \Re^n$ be a continuous mapping with $b > a > -\infty$ and $r > 0$. By $x_t$ we denote the "$r$-history" of $x$ at time $t \in [a,b)$, i.e., $(x_t)(\theta) := x(t+\theta); \theta \in [-r,0]$. Notice that $x_t \in C^0([-r,0];\Re^n)$.

* By $K$ we denote the set of increasing, continuous functions $\rho: \Re_+ \to \Re_+$ with $\rho(0) = 0$. We say that a function $\rho \in K$ is of class $K_\infty$ if $\lim_{s \to +\infty} \rho(s) = +\infty$.

* Let $D \subseteq \Re^l$ be a non-empty set and $I \subseteq \Re_+$ an interval. By $L^\infty_{loc}(I;D)$ we denote the class of Lebesgue measurable and locally bounded mappings $d: \Re_+ \to D$. Notice that by $\sup_{\tau \in [0,t]}(|d(\tau)|)$ we do not mean the essential supremum of $d: \Re_+ \to D$ on $[0,t]$ but the actual supremum of $d: \Re_+ \to D$ on $[0,t]$.

## 2. Assumptions and Main Result

In the present work we study systems of the form (1.1) under the following assumptions:

**(H1)** *The mappings $f: C^0([-r,0]; \Re^n) \times U \times D \to \Re^n$, $h: C^0([-r,0]; \Re^n) \to \Re^k$ are continuous and satisfy the following properties: (i) $f(0,0,0) = 0$, $h(0) = 0$, (ii) for every bounded $\Omega \subset C^0([-r,0]; \Re^n) \times U \times D$ the image set $f(\Omega) \subset \Re^n$ is bounded, (iii) for every bounded $\Omega \subset C^0([-r,0]; \Re^n)$ the image set $h(\Omega) \subset \Re^k$ is bounded, and (iv) for every bounded $S \subset C^0([-r,0]; \Re^n) \times U \times D$, there exists a constant $L_S \geq 0$ such that*

$$(x(0) - \overline{x}(0))^T (f(x,u,d) - f(\overline{x},u,d)) \leq L_S \|x - \overline{x}\|^2$$

$$\forall (x,u,d) \in S, \forall (\overline{x},u,d) \in S$$

**(H2)** *System (1.1) is Forward Complete, i.e., for every $x_0 \in C^0([-r,0]; \Re^n)$ and for every $u \in L^\infty_{loc}(\Re_+; U)$, $d \in L^\infty_{loc}(\Re_+; D)$ the solution of (1.1) with initial condition $x_0$, corresponding to inputs $u, d$, exists for all $t \geq 0$.*

Assumption (H1) is a standard assumption for time-delay systems that guarantees existence and uniqueness of solutions for system (1.1), i.e. guarantees that for every $x_0 \in C^0([-r,0]; \Re^n)$ and for every $u \in L^\infty_{loc}(\Re_+; U)$, $d \in L^\infty_{loc}(\Re_+; D)$ there exists $t_{\max} \in (0, +\infty]$ and a unique continuous mapping $x: [-r, t_{\max}) \to \Re^n$ which is absolutely continuous on $[0, t_{\max})$ and satisfies $x(\theta) = x_0(\theta)$ for $\theta \in [-r, 0]$ and (1.1) for $t \in [0, t_{\max})$ a.e.. This mapping $x: [-r, t_{\max}) \to \Re^n$ is the solution of (1.1) with initial condition $x_0$, corresponding to inputs $u, d$. Assumption (H2) guarantees that $t_{\max} = +\infty$ for every $x_0 \in C^0([-r,0]; \Re^n)$ and for every $u \in L^\infty_{loc}(\Re_+; U)$, $d \in L^\infty_{loc}(\Re_+; D)$.

The following assumption plays a crucial role in what follows.

**(H3)** *There exists a continuous mapping $R: C^0([-r,0]; \Re^n) \times U \times D \to \Re^k$ with the following property: for every $x_0 \in C^0([-r,0]; \Re^n)$, $u \in L^\infty_{loc}(\Re_+; U)$, $d \in L^\infty_{loc}(\Re_+; D)$ the unique solution $x$ of (1.1) with initial condition $x_0$, corresponding to inputs $u, d$, satisfies for $t \geq 0$ a.e. the following equation:*

$$\frac{d}{dt}(h(x_t)) = R(x_t, u(t), d(t)) \qquad (2.1)$$

*Moreover, there exists a constant $L \geq 0$ and a function $\kappa \in K$ such that the following inequality holds for all $x, \overline{x} \in C^0([-r,0]; \Re^n)$, $(u,d) \in U \times D$:*

$$|R(\overline{x}, u, 0) - R(x, u, d)| \leq L\|\overline{x} - x\| + \kappa(|d|) \qquad (2.2)$$

Assumption (H3) essentially requires that the derivative of the output of system (1.1) exists and is expressed by the globally Lipschitz (with respect to $x$) mapping $R: C^0([-r,0]; \Re^n) \times U \times D \to \Re^k$. Clearly, not every nonlinear time-delay system satisfies (H3). Nevertheless, the class of systems satisfying Assumption (H3) is very wide, including many systems of practical interest e.g. (1.6).

The notion of the Robust Exponential Observer (REO) for system (1.1) is crucial to the development of the main results of the present work and it is given in the following definition.

**Definition 2.1 (Robust Exponential Observer):** *Consider the following system*

$$\dot{z} = F(z_t, y, u), z \in \Re^l$$
$$\hat{x}_t = \Psi(z_t), \hat{x} \in \Re^n \tag{2.3}$$

*where* $F: C^0([-r,0];\Re^l) \times \Re^k \times U \to \Re^l$, $\Psi: C^0([-r,0];\Re^l) \to C^0([-r,0];\Re^n)$ *are continuous mappings with* $F(0,0,0) = 0$, $\Psi(0) = 0$. *Suppose that the mapping* $F: C^0([-r,0];\Re^l) \times \Re^k \times U \to \Re^l$ *is such that, for every bounded* $\Omega \subset C^0([-r,0];\Re^l) \times \Re^k \times U$ *the image set* $F(\Omega) \subset \Re^l$ *is bounded and there exists a constant* $L_\Omega \geq 0$ *such that*

$$(z(0) - \bar{z}(0))^T (F(z,y,u) - F(\bar{z},y,u)) \leq L_\Omega \|z - \bar{z}\|^2$$
$$\forall (z, y, u) \in \Omega, \forall (\bar{z}, y, u) \in \Omega$$

*System (2.3) is called a **Robust Exponential Observer (REO)** for system (1.1), if there exist constants* $\gamma, \sigma > 0$ *and function* $g \in K$, $a \in K_\infty$ *such that for every* $(x_0, z_0) \in C^0([-r,0];\Re^n) \times C^0([-r,0];\Re^l)$ *and* $u \in L^\infty_{loc}(\Re_+;U)$, $d \in L^\infty_{loc}(\Re_+;D)$, $\xi \in L^\infty_{loc}(\Re_+;\Re^k)$, *the solution* $(x(t), z(t))$ *of*

$$\dot{x} = f(x_t, u, d)$$
$$\dot{z} = F(z_t, h(x_t) + \xi, u) \tag{2.4}$$
$$\hat{x}_t = \Psi(z_t)$$

*with initial condition* $(x_0, z_0)$, *corresponding to inputs* $u, d, \xi$, *exists for all* $t \geq 0$ *and satisfies the following estimate:*

$$\|\hat{x}_t - x_t\| \leq \exp(-\sigma t) a(\|x_0\| + \|z_0\|) + \gamma \sup_{0 \leq s \leq t} (|\xi(s)| \exp(-\sigma(t-s))) + \sup_{0 \leq s \leq t} (g(|d(s)|)), \forall t \geq 0 \tag{2.5}$$

At this point, it should be noticed that the way that the inputs $d$ and $\xi$ enter the Input-to-Output (IOS) Stability estimate (2.5) is different. While the input $d$ comes in estimate (2.5) through a (possibly) nonlinear gain function, the input $\xi$ appears in estimate (2.5) with a linear gain and with a fading memory effect (see [24]). This difference is important and allows, in what follows, the construction of sampled-data observers.

Besides the fact that Definition 2.1 introduces the notion of REO for systems with state delays, there are important differences between the notion of a REO in Definition 2.1 and other similar notions that were used in the literature (for systems described by ODEs; see [4,23]):
1) In Definition 2.1, the effect of disturbances is explicitly taken into account (see the term $\sup_{0 \leq s \leq t} g(|d(s)|)$ in estimate (2.5)), while in other similar notions in the literature no disturbances are assumed to act on the system.
2) In Definition 2.1, the IOS estimate (2.5) is assumed to hold uniformly for all inputs $u \in L^\infty_{loc}(\Re_+;U)$, while in other similar notions in the literature either there is no control input $u$ or the sup-norm of $u$ appears in the corresponding observer error estimate. This difference is important when the observer is to be used in conjunction with a state feedback control law for the dynamic stabilization of the system.

We are now in a position to state our main result.

**Theorem 2.2 (Sampled-Data Observer Design):** *Consider system (1.1) under Assumptions (H1), (H2), (H3) and suppose that system (2.3) is a REO for system (1.1). Moreover, suppose that for every bounded $S \subset C^0([-r,0];\Re^n) \times U$, there exists a constant $L_S \geq 0$ such that*

$$\begin{aligned} |R(\Psi(z),u,0) - R(\Psi(\bar{z}),u,0)| &\leq L_S \|z - \bar{z}\| \\ \forall (z,u) \in S, \ \forall (\bar{z},u) \in S \end{aligned} \quad (2.6)$$

*Let $\delta > 0$ and $\omega \in (0, \sigma]$ be constants that satisfy*

$$\delta < \frac{1}{\omega} \ln\left(1 + \frac{\omega}{\gamma L}\right) \quad (2.7)$$

*Then for every sampling sequence $\{\tau_i\}_{i=0}^{\infty}$ with $\tau_0 = 0$, $\lim(\tau_i) = +\infty$, $0 < \tau_{i+1} - \tau_i \leq \delta$ for $i = 0, 1, \ldots,$ for every $(x_0, z_0) \in C^0([-r,0];\Re^n) \times C^0([-r,0];\Re^l)$ and $u \in L^{\infty}_{loc}(\Re_+;U)$, $d \in L^{\infty}_{loc}(\Re_+;D)$, $\xi \in L^{\infty}_{loc}(\Re_+;\Re^k)$, the solution $(x(t), z(t), w(t))$ of (1.1) with*

$$\begin{aligned} \dot{z}(t) &= F(z_t, w(t), u(t)), \ t \geq 0 \\ \dot{w}(t) &= R(\Psi(z_t), u(t), 0), \ t \in [\tau_i, \tau_{i+1}) \\ \hat{x}_t &= \Psi(z_t), \ t \geq 0 \end{aligned} \quad (2.8)$$

$$w(\tau_i) = h(x_{\tau_i}) + \xi(\tau_i) \quad (2.9)$$

*initial condition $(x_0, z_0)$, corresponding to inputs $u, d, \xi$, exists for all $t \geq 0$ and satisfies the following estimate*

$$\begin{aligned} \|\hat{x}_t - x_t\| &\leq (1-B)^{-1} \exp(-\omega t) a\left(\|x_0\| + \|z_0\|\right) \\ &\quad + (1-B)^{-1} \gamma \exp(\omega \delta) \sup_{0 \leq s \leq t}\left(|\xi(s)| \exp(-\omega(t-s))\right) + \sup_{0 \leq s \leq t}\left(\tilde{g}(|d(s)|)\right) \end{aligned} \quad (2.10)$$

*where $\tilde{g}(s) := (1-B)^{-1}(g(s) + \gamma \delta \kappa(s))$ and $B := \gamma L \frac{\exp(\omega \delta) - 1}{\omega} < 1$.*

**Remark 2.3:** (a) The sampled-data observer (2.8), (2.9) is simply the REO (2.3) with the unavailable output signal replaced by the signal produced by the inter-sample output predictor

$$\begin{aligned} \dot{w}(t) &= R(\Psi(z_t), u(t), 0), \ t \in [\tau_i, \tau_{i+1}) \\ w(\tau_i) &= h(x_{\tau_i}) + \xi(\tau_i) \end{aligned}$$

(b) Notice that (2.10) guarantees the IOS property for the output map $Y = \hat{x}_t - x_t$ from the inputs $d, \xi$, i.e. from the inputs expressing the effect of modeling errors and measurement noise, respectively. However, a comparison of (2.5) and (2.10) shows that the input gains are higher for the sampled-data observer (2.8), (2.9) than the continuous-time REO (2.3). It is clear that sampling makes the observer more sensitive to modeling errors and measurement noise.

**(c)** The sampled-data observer (2.8), (2.9) is a hybrid observer with delays. For each sampling sequence $\{\tau_i\}_{i=0}^{\infty}$ with $\tau_0 = 0$, $\lim(\tau_i) = +\infty$, $0 < \tau_{i+1} - \tau_i \leq \delta$ for $i = 0, 1, \ldots$, for each $(x_0, z_0) \in C^0([-r,0]; \Re^n) \times C^0([-r,0]; \Re^l)$ and $u \in L^{\infty}_{loc}(\Re_+; U)$, $d \in L^{\infty}_{loc}(\Re_+; D)$, $\xi \in L^{\infty}_{loc}(\Re_+; \Re^k)$, the solution $(x(t), z(t), w(t))$ of (1.1) with (2.8), (2.9) with initial condition $(x_0, z_0)$, corresponding to inputs $u, d, \xi$, is produced by the following algorithm:

Step $i \geq 0$:
1) Given $x_{\tau_i}$, calculate $x_t$ for $t \in (\tau_i, \tau_{i+1}]$ from (1.1) and calculate $w(\tau_i) = h(x_{\tau_i}) + \xi(\tau_i)$,
2) Given $z_{\tau_i}$ and $w(\tau_i)$, calculate $z_t$ for $t \in (\tau_i, \tau_{i+1}]$ and $w(t)$ for $t \in (\tau_i, \tau_{i+1})$ as the solution of the time-delay system $\dot{z}(t) = F(z_t, w(t), u(t))$ and $\dot{w}(t) = R(\Psi(z_t), u(t), 0)$,
3) Compute the output trajectory $\hat{x}_t$, for $t \in (\tau_i, \tau_{i+1}]$ using the equation $\hat{x}_t = \Psi(z_t)$

**(d)** Despite the hybrid nature of the observer (2.8)-(2;9), the trajectory of the estimated state features continuity. The proof of Theorem 2.2 is based on a small-gain argument (see Section 5). It is therefore expected that the observer error estimate (2.10) and the upper bound for the diameter of the sampling sequence $\delta > 0$ given by (2.7) are conservative. However, formulas (2.7), (2.10) are useful because they indicate which parameters affect the performance of the observer and (qualitatively) how the upper bound for the diameter of the sampling sequence depends on the parameters of the system.

**(e)** Since the mapping $(0, \sigma] \ni \omega \to \frac{1}{\omega} \ln\left(1 + \frac{\omega}{\gamma L}\right)$ is decreasing with $\lim_{\omega \to 0^+}\left(\frac{1}{\omega} \ln\left(1 + \frac{\omega}{\gamma L}\right)\right) = \frac{1}{\gamma L}$, it is clear from (2.7) that: (i) Theorem 2.2 requires sampling sequences with diameter $\delta > 0$ being less than $\frac{1}{\gamma L}$, and (ii) the smaller the diameter $\delta > 0$ of the sampling sequence is, the larger the constant $\omega > 0$ is, i.e., convergence is faster for a smaller diameter $\delta > 0$ of the sampling sequence in the absence of modeling errors and measurement noise (recall (2.10)).

**(f)** In general, the constants $\gamma$ and $L$ depend on the value of the maximum delay $r$. Therefore, inequality (2.7) provides a useful relation between the diameter of the sampling sequence $\delta$ and the delay $r$.

For the design of observed-based, output feedback we need a stabilizability assumption.

**(H4)** *The equalities $f(0,0,0) = 0$ and $U = \Re^m$ hold. Moreover there exist a function $\bar{\kappa} \in K$, constants $\sigma, M > 0$ and a functional $\tilde{k}: C^0([-r,0]; \Re^n) \to \Re^m$ with $\tilde{k}(0) = 0$ and a constant $\bar{L} > 0$ such that the inequalities $\left|\tilde{k}(\bar{x}) - \tilde{k}(x)\right| + \left|f(\bar{x}, \bar{u}, 0) - f(x, u, 0)\right| \leq \bar{L}\|\bar{x} - x\| + \bar{L}|\bar{u} - u|$, $|f(x,u,d) - f(x,u,0)| \leq \bar{\kappa}(|d|)$ hold for all $x, \bar{x} \in C^0([-r,0]; \Re^n)$, $\bar{u}, u \in \Re^m$, $d \in D$ and such that for every $x_0 \in C^0([-r,0]; \Re^n)$, the solution $x(t)$ of*

$$\dot{x}(t) = f(x_t, u(t), 0)$$
$$u(t) = \tilde{k}(x_t) \tag{2.11}$$

*with initial condition $x_0$ exists for all $t \geq 0$ and satisfies the following estimate*

$$\|x_t\| \leq M \exp(-\sigma t)\|x_0\|, \quad \forall t \geq 0 \tag{2.12}$$

When Assumption (H4) holds then we obtain the following stabilization result.

**Corollary 2.4 (Global Stabilization by Means of Observer-Based Sampled-Data Output Feedback):** *Consider system (1.1) under Assumptions (H1), (H2), (H3), (H4) and suppose that system (2.3) is a REO for system (1.1). Moreover, suppose that for every bounded $S \subset C^0([-r,0];\Re^n) \times U$, there exists a constant $L_S \geq 0$ such that (2.6) holds. Then there exist constants $\bar{\delta}, \omega, \hat{\gamma} > 0$ and functions $\hat{g}, \hat{a} \in K$ such that for every sampling sequence $\{\tau_i\}_{i=0}^{\infty}$ with $\tau_0 = 0$, $\lim(\tau_i) = +\infty$, $0 < \tau_{i+1} - \tau_i \leq \bar{\delta}$ for $i = 0,1,...$, for every $(x_0, z_0) \in C^0([-r,0];\Re^n) \times C^0([-r,0];\Re^l)$ and $\xi \in L^{\infty}_{loc}(\Re_+;\Re^k)$, $d \in L^{\infty}_{loc}(\Re_+;D)$, the solution $x(t)$ of (1.1) with (2.8), (2.9) and*

$$u(t) = \tilde{k}(\hat{x}_{\tau_i}), t \in [\tau_i, \tau_{i+1}) \tag{2.13}$$

*initial condition $(x_0, z_0)$, corresponding to inputs $d, \xi$, exists for all $t \geq 0$ and satisfies the following estimate*

$$\|x_t\| + \|\hat{x}_t\| \leq \exp(-\omega t)\hat{a}(\|x_0\| + \|z_0\|) + \hat{\gamma} \sup_{0 \leq s \leq t}(|\xi(s)|) + \sup_{0 \leq s \leq t}(\hat{g}(|d(s)|)), \forall t \geq 0 \tag{2.14}$$

*Moreover, if $a \in K_{\infty}$ (the function involved in (2.5)) is linear then $\hat{a}$ is linear too.*

## 3. Hyperbolic PDE-ODE Interconnections

*3.A. The General Case*

When a plant is interconnected with a transport process then we can obtain a system of the form (1.1) with distributed delays. This is the reason that in this section we consider initial-boundary value problems of the form (1.2), (1.3), (1.4) with initial condition

$$v[0] = v_0, \bar{x}(0) = \bar{x}_0 \tag{3.1}$$

where $c > 0$ is a constant, $\bar{x}(t) \in \Re^{\bar{n}}$, $v(t,z) \in \Re$ are the states, $u \in C^0(\Re_+;\Re^m)$ is an external input, $a \in C^0([0,1])$, $b_i \in C^1([0,1])$ ($i=1,...,N_1$), $\gamma_i \in C^0([0,1] \times [0,1])$, $\beta_i \in C^1([0,1])$ ($i=1,...,N_2$), $\bar{x}_0 \in \Re^n$, $v_0 \in C^1([0,1])$ satisfy the compatibility conditions $cv_0'(0) = g(0,\bar{x}_0) + \sum_{i=1}^{N_1} b_i(0)v_0(z_i) + \sum_{i=1}^{N_2} \beta_i(0)\int_0^1 \gamma_i(s)v_0(s)ds$, $v_0(0) = 0$ and the following assumption holds for the mappings $\tilde{F}:\Re^{\bar{n}} \times C^0([0,1]) \times \Re^m \to \Re^{\bar{n}}$, $g:[0,1] \times \Re^{\bar{n}} \to \Re$:

**(A1)** $\tilde{F}:\Re^{\bar{n}} \times C^0([0,1]) \times \Re^m \to \Re^{\bar{n}}$, $g \in C^1([0,1] \times \Re^{\bar{n}};\Re)$ *are continuous mappings with* $\tilde{F}(0,0,0) = 0$, $g(z,0) = 0$ *for all $z \in [0,1]$, for which there exists a constant $L > 0$ such that the inequalities $|\tilde{F}(x,v,u) - \tilde{F}(y,w,u)| \leq L\|v-w\|_{\infty} + L|x-y|$, $\max_{0 \leq z \leq 1}(|g(z,x) - g(z,y)|) \leq L|x-y|$ hold for all $v, w \in C^0([0,1])$, $x, y \in \Re^{\bar{n}}$, $u \in \Re^m$.*

Under Assumption (A1), Theorem 2.2 on page 22 in [28] guarantees that for every $v_0 \in C^1([0,1])$, $\bar{x}_0 \in \Re^{\bar{n}}$ and $u \in C^0(\Re_+;\Re^m)$ with $cv_0'(0) = g(0,\bar{x}_0) + \sum_{i=1}^{N_1} b_i(0)v_0(z_i) + \sum_{i=1}^{N_2} \beta_i(0)\int_0^1 \gamma_i(s)v_0(s)ds$, $v_0(0) = 0$, there exist unique mappings $v \in C^1(\Re_+ \times [0,1])$ and $\bar{x} \in C^1(\Re_+; \Re^{\bar{n}})$ satisfying (1.2), (1.3), (1.4), (3.1). Moreover, the solution satisfies the formula

$$\begin{aligned}v(t,z) &= v_0(\max(0, z-ct))\exp\left(c^{-1}\int_{\max(0,z-ct)}^{z} a(w)dw\right) \\ &+ c^{-1}\int_{\max(0,z-ct)}^{z} \exp\left(c^{-1}\int_{l}^{z} a(w)dw\right) g(l, \bar{x}(t-c^{-1}z+c^{-1}l))dl \\ &+ c^{-1}\sum_{i=1}^{N_1}\int_{\max(0,z-ct)}^{z} \exp\left(c^{-1}\int_{l}^{z} a(w)dw\right) b_i(l) v(t-c^{-1}z+c^{-1}l, z_i)dl \\ &+ c^{-1}\sum_{i=1}^{N_2}\int_{\max(0,z-ct)}^{z} \exp\left(c^{-1}\int_{l}^{z} a(w)dw\right) \beta_i(l) \int_0^1 \gamma_i(s) v(t-c^{-1}z+c^{-1}l, s)ds\,dl \\ &\text{for } (t,z) \in \Re_+ \times [0,1]\end{aligned} \quad (3.2)$$

For such systems, the output $y(t) = (y_1(t), ..., y_k(t))^T \in \Re^k$ is given by (1.5), where $q_j \in \Re^{\bar{n}}$ ($j = 1,...,k$) are constant vectors, $\bar{b}_{j,i} \in \Re$ ($j = 1,...,k$, $i = 1,...,N_1$) and $\bar{\beta}_{j,i} \in \Re$ ($j = 1,...,k$, $i = 1,...,N_2$) are constants.

Therefore, for $t \geq r$, where

$$r = c^{-1} \quad (3.3)$$

we obtain from (3.2), (3.3) and (1.5) the delay system (1.2) with

$$v(t,z) = \int_{t-rz}^{t} \frac{C(z)}{C(z-c(t-p))} \bar{g}(z-c(t-p), \bar{x}(p), \eta(p), \lambda(p))dp$$
$$\text{for } t \geq r, \; z \in [0,1] \quad (3.4)$$

with

$$C(z) := \exp\left(r\int_0^z a(s)ds\right) \quad (3.5)$$

$$\bar{g}(z, \bar{x}, \eta, \lambda) := g(z, \bar{x}) + \sum_{i=1}^{N_1} b_i(z)\eta_i + \sum_{i=1}^{N_2} \beta_i(z)\lambda_i \quad (3.6)$$

$$\eta_i(t) = v(t, z_i) = \int_{t-rz_i}^{t} \frac{C(z_i)}{C(z_i - c(t-p))} \bar{g}(z_i - c(t-p), \bar{x}(p), \eta(p), \lambda(p))dp$$
$$\text{for } t \geq r, \; i = 1, ..., N_1 \quad (3.7)$$

$$\lambda_i(t) = \int_0^1 \gamma_i(z) v(t, z)dz = \int_0^1 \int_{t-rz}^{t} \frac{\gamma_i(z) C(z)}{C(z - c(t-p))} \bar{g}(z - c(t-p), \bar{x}(p), \eta(p), \lambda(p))dp\,dz$$
$$\text{for } t \geq r, \; i = 1, ..., N_2 \quad (3.8)$$

and output given by

$$y_j(t) = q_j^T \bar{x}(t) + \sum_{i=1}^{N_1} \bar{b}_{j,i} \eta_i(t) + \sum_{i=1}^{N_2} \bar{\beta}_{j,i} \lambda_i(t), \text{ for } t \geq r, \; j = 1,...,k \quad (3.9)$$

At this point, we need the following technical assumption:

**(A2)** *There exists a constant $G > 0$ such that for every $v_0 \in C^1([0,1])$, $\bar{x}_0 \in \Re^{\bar{n}}$ with*

$$cv_0'(0) = g(0, \bar{x}_0) + \sum_{i=1}^{N_1} b_i(0) v_0(z_i) + \sum_{i=1}^{N_2} \beta_i(0) \int_0^1 \gamma_i(s) v_0(s) ds, \qquad v_0(0) = 0, \qquad \text{there exist}$$

$\bar{x} \in C^0([-r,0]; \Re^{\bar{n}})$, $\eta \in C^0([-r,0]; \Re^{N_1})$, $\lambda \in C^0([-r,0]; \Re^{N_2})$ with $\bar{x}(0) = \bar{x}_0$, $\eta_i(0) = v_0(z_i)$

($i = 1, ..., N_1$), $\lambda_i(0) = \int_0^1 \gamma_i(z) v_0(z) dz$ ($i = 1, ..., N_2$) *and* $\|\bar{x}\| + \|\eta\| + \|\lambda\| \leq G(\|v_0\|_\infty + \|v_0'\|_\infty + |\bar{x}_0|)$ *that satisfy*

$$v_0(s) = \int_{-rs}^{0} \frac{C(s)}{C(cp+s)} \tilde{g}(cp+s, \bar{x}(p), \eta(p), \lambda(p)) dp, \text{ for all } s \in [0,1] \qquad (3.10)$$

If Assumption (A2) holds then equations (3.4), (3.7), (3.8), (3.9) are valid for all $(t,z) \in \Re_+ \times [0,1]$. Moreover, in this case and if the output is sampled, then Theorem 2.2 can be used for sampled-data observer design. We define

$$x = \begin{pmatrix} \bar{x} \\ \eta \\ \lambda \end{pmatrix} \in \Re^n, \text{ with } n = \bar{n} + N_1 + N_2, \qquad (3.11)$$

$$f(x,u,d) := \begin{pmatrix} f_1(x,u,d) \\ f_2(x) \\ f_3(x) \end{pmatrix}, \quad f_2(x) := \begin{pmatrix} f_{2,1}(x) \\ \vdots \\ f_{2,N_1}(x) \end{pmatrix}, \quad f_3(x) := \begin{pmatrix} f_{3,1}(x) \\ \vdots \\ f_{3,N_2}(x) \end{pmatrix} \qquad (3.12)$$

for all $x \in C^0([-r,0]; \Re^n)$, $u \in \Re^m$ and $d \in \Re^q$ by means of the following equations:

$$f_1(x,u,d) := \tilde{F}(\bar{x}(0), v, u) + d, \text{ with } v(z) = \int_{-rz}^{0} \frac{C(z)}{C(cp+z)} \tilde{g}(cp+z, x(p)) dp, \; z \in [0,1] \quad (3.13)$$

$$f_{2,i}(x) := \bar{g}(z_i, x(0)) - C(z_i) \bar{g}(0, x(-rz_i)) + \int_{-rz_i}^{0} \frac{C(z_i)}{C(z_i + cp)} \tilde{g}(z_i + cp, x(p)) dp, \text{ for } i = 1, ..., N_1 \quad (3.14)$$

$$f_{3,i}(x) := \int_0^1 \gamma_i(z) (\bar{g}(z, x(0)) - C(z) \bar{g}(0, x(-rz))) dz + \int_0^1 \int_{-rz}^{0} \frac{\gamma_i(z) C(z)}{C(z+cp)} \tilde{g}(z+cp, x(p)) dp dz$$

for $i = 1, ..., N_2$ \hfill (3.15)

$$h_j(x) := q_j^T \bar{x}(0) + \sum_{i=1}^{N_1} \bar{b}_{j,i} \eta_i(0) + \sum_{i=1}^{N_2} \bar{\beta}_{j,i} \lambda_i(0), \text{ for } j = 1, ..., k \qquad (3.16)$$

where

$$\tilde{g}(z, \bar{x}, \eta, \lambda) := a(z) g(z, \bar{x}) - c \frac{\partial g}{\partial z}(z, \bar{x}) + \sum_{j=1}^{N_1} \tilde{b}_j(z_i) \eta_j + \sum_{j=1}^{N_2} \tilde{\beta}_j(z_i) \lambda_j \qquad (3.17)$$

$$\tilde{b}_j(z) := a(z) b_j(z) - c b_j'(z), \; i = 1, ..., N_1 \qquad (3.18)$$

$$\tilde{\beta}_j(z) := a(z) \beta_j(z) - c \beta_j'(z), \; i = 1, ..., N_2 \qquad (3.19)$$

The proof of the following lemma is almost trivial and is omitted.

**Lemma 3.1:** *Suppose that Assumptions (A1), (A2) hold for system (1.2), (1.3), (1.4). Moreover, suppose that the following assumption holds:*

**(A3)** *There exists a constant $L>0$ such that the inequality $\max_{0\leq z\leq 1}\left(\left|\frac{\partial g}{\partial z}(z,x)-\frac{\partial g}{\partial z}(z,y)\right|\right)\leq L|x-y|$ holds for all $x,y\in\Re^n$.*

*Consider system (1.1) with $D=\{0\}$, where $f:C^0([-r,0];\Re^n)\times\Re^m\times D\to\Re^n$, $h:C^0([-r,0];\Re^n)\to\Re^k$ are defined by (3.12), (3.13), (3.14), (3.15), (3.16). Then Assumptions (H1), (H2), (H3) for system (1.1) hold with $U=\Re^m$ and $R:C^0([-r,0];\Re^n)\times\Re^m\times D\to\Re^k$ defined for all $(x,u,d)\in C^0([-r,0];\Re^n)\times\Re^m\times D$ by*

$$R(x,u,d)=\begin{pmatrix} q_1^T f_1(x,u,d)+\sum_{i=1}^{N_1}\bar{b}_{1,i}f_{2,i}(x)+\sum_{i=1}^{N_2}\bar{\beta}_{1,i}f_{3,i}(x) \\ \vdots \\ q_k^T f_1(x,u,d)+\sum_{i=1}^{N_1}\bar{b}_{k,i}f_{2,i}(x)+\sum_{i=1}^{N_2}\bar{\beta}_{k,i}f_{3,i}(x) \end{pmatrix} \quad (3.20)$$

Based on Theorem 2.2 and Lemma 3.1, we are in a position to show the following result.

**Theorem 3.2 (Sampled-Data Observer for Hyperbolic PDE-ODE Loops):** *Suppose that Assumptions (A1), (A2), (A3) hold for system (1.2), (1.3), (1.4). Consider system (1.1) with $D=\{0\}$, where $f:C^0([-r,0];\Re^n)\times\Re^m\times D\to\Re^n$, $h:C^0([-r,0];\Re^n)\to\Re^k$ are defined by (3.12), (3.13), (3.14), (3.15), (3.16). Suppose that system (2.3) is a REO for system (1.1). Moreover, suppose that for every bounded $S\subset C^0([-r,0];\Re^n)\times\Re^m$, there exists a constant $L_S\geq 0$ such that (2.6) holds with $R:C^0([-r,0];\Re^n)\times\Re^m\times D\to\Re^k$ being defined by (3.20). Then there exist constants $\delta,P,\omega>0$ and a function $\bar{a}\in K_\infty$ such that for every sampling sequence $\{\tau_i\}_{i=0}^\infty$ with $\tau_0=0$, $\lim(\tau_i)=+\infty$, $0<\tau_{i+1}-\tau_i\leq\delta$ for $i=0,1,...$, for every $v_0\in C^1([0,1])$, $\bar{x}_0\in\Re^{\bar{n}}$ with $cv_0'(0)=g(0,\bar{x}_0)+\sum_{i=1}^{N_1}b_i(0)v_0(z_i)+\sum_{i=1}^{N_2}\beta_i(0)\int_0^1\gamma_i(s)v_0(s)ds$, $v_0(0)=0$, $z_0\in C^0([-r,0];\Re^l)$ and $u\in C^0(\Re_+;\Re^m)$, $\xi\in L_{loc}^\infty(\Re_+;\Re^k)$, the solution $(x(t),v[t],z(t),w(t))$ of (1.2), (1.3), (1.4) together with (2.8) and*

$$w_j(\tau_i)=\xi(\tau_i)+q_j^T\bar{x}(\tau_i)+\sum_{l=1}^{N_1}\bar{b}_{j,l}v(\tau_i,z_l)+\sum_{l=1}^{N_2}\bar{\beta}_{j,l}\int_0^1\gamma_l(s)v(\tau_i,s)ds,\ j=1,...,k \quad (3.21)$$

$$\hat{v}(t,z)=\int_{t-rz}^t\frac{C(z)}{C(z-c(t-p))}\bar{g}(z-c(t-p),\hat{x}(p))dp,\text{ for } t\geq 0,\ z\in[0,1] \quad (3.22)$$

*initial condition $(\bar{x}_0,v_0,z_0)$, corresponding to inputs $u,\xi$, exists for all $t\geq 0$ and satisfies the following estimate:*

$$\begin{aligned}&|\bar{x}(t)-H\hat{x}(t)|+\|v[t]-\hat{v}[t]\|_\infty\\&\leq\exp(-\omega t)\bar{a}\left(\|v_0\|_\infty+\|v_0'\|_\infty+|\bar{x}_0|+\|z_0\|\right)+P\sup_{0\leq s\leq t}\left(|\xi(s)|\exp(-\omega(t-s))\right)\end{aligned} \quad (3.23)$$

where $H \in \Re^{\bar{n} \times n}$ is the matrix for which the relation

$$\bar{x} = H \begin{pmatrix} \bar{x} \\ \eta \\ \lambda \end{pmatrix} \quad (3.24)$$

holds for all $\bar{x} \in \Re^{\bar{n}}$, $\eta \in \Re^{N_1}$, $\lambda \in \Re^{N_2}$. Moreover, if $a \in K_\infty$ (the function involved in (2.5)) is linear then $\bar{a}$ is linear.

Theorem 3.2 can be used as a tool for the construction of sampled-data observers for hyperbolic PDE-ODE loops. It allows the transformation of the initial sampled-data observer design problem for a PDE-ODE loop to the construction of a REO for a specific delay system.

**Remark:** It should be noted that under Assumptions (A1), (A2), (A3), system (1.2), (1.3), (1.4) is *not* equivalent to the delay system

$$\frac{d\bar{x}}{dt}(t) = f_1(\bar{x}_t, \eta_t, \lambda_t, u(t), d(t))$$
$$\frac{d\eta}{dt}(t) = f_2(\bar{x}_t, \eta_t, \lambda_t) \quad (3.25)$$
$$\frac{d\lambda}{dt}(t) = f_3(\bar{x}_t, \eta_t, \lambda_t)$$

Indeed, for every solution of system (1.2), (1.3), (1.4) there exists a solution of system (3.25) for which (3.4), (3.7), (3.8) hold. However, not every solution of system (3.25) provides a solution of system (1.2), (1.3), (1.4) by means of (3.4), (3.7), (3.8). For such a thing to happen, the initial condition of the solution $(\bar{x}_t, \eta_t, \lambda_t)$ of system (3.25) has to satisfy the additional equations

$$\eta_i(0) = \int_{-rz_i}^{0} \frac{C(z_i)}{C(z_i + cp)} \bar{g}(z_i + cp, \bar{x}(p), \eta(p), \lambda(p)) dp, \quad \text{for} \quad i = 1, \ldots, N_1 \quad \text{and}$$

$$\lambda_i(0) = \int_0^1 \int_{-rz}^0 \frac{\gamma_i(z) C(z)}{C(z + cp)} \bar{g}(z + cp, \bar{x}(p), \eta(p), \lambda(p)) dp dz, \quad \text{for} \quad i = 1, \ldots, N_2.$$ In other words, system (1.2), (1.3), (1.4) is *immersed into* system (3.25).

*3.B. Application to a Chemical Reactor*

The model of a chemical reactor with an exothermic chemical reaction taking place in it and a cooling jacket with negligible axial heat conduction of the cooling medium is given in Chapter 2 and Chapter 8 of [28]:

$$\frac{d\bar{x}}{dt}(t) = \theta(\bar{x}(t)) - (\mu + 1)\bar{x}(t) + \mu \int_0^1 v(t, z) dz \quad (3.26)$$

$$\frac{\partial v}{\partial t}(t, z) + c \frac{\partial v}{\partial z}(t, z) = -\zeta v(t, z) + \zeta \bar{x}(t) \quad (3.27)$$

with (1.4), where $\mu, \zeta > 0$ are constants, $\bar{x}(t)$ is the (dimensionless) reactor outlet temperature, $v(t, z)$ is the (dimensionless) temperature of the cooling medium at position $z \in [0,1]$ in the jacket ($z = 0$ is the entrance of the jacket and $z = 1$ is the exit of the jacket) and $\theta : \Re \to \Re$ is a $C^1$ globally Lipschitz function with $\theta(0) = 0$. When the reactor inlet temperature is not constant then equation (3.26) is modified as follows

$$\frac{d\bar{x}}{dt}(t) = \theta(\bar{x}(t)) - (\mu+1)\bar{x}(t) + \mu\int_0^1 v(t,z)dz + u(t) \tag{3.28}$$

where $u(t) \in \Re$ is the input that expresses the variation of the reactor inlet temperature. Moreover, it is indeed the case that the measured temperature is the temperature of the cooling medium at the exit of the jacket, i.e., the measured output is $y(t) = v(t,1)$. The reactor model (1.4), (3.27), (3.28), is a system of the form (1.2), (1.3), (1.4), (1.5) with $a(z) \equiv -\zeta$, $\bar{n} = 1$, $N_1 = 1$, $z_1 = 1$, $k = 1$, $b_1(z) \equiv 0$, $g(z,\bar{x}) = \zeta\bar{x}$, $\bar{b}_{1,1} = 1$, $q_1 = 0$, $\tilde{F}(\bar{x},v,u) = \theta(\bar{x}) - (\mu+1)\bar{x} + \mu\int_0^1 v(z)dz + u$ ($N_2$, $\gamma_i$, $\beta_i$, $\bar{\beta}_{1,i}$ are irrelevant).

Assumptions (A1), (A3) are automatically verified. Moreover, Assumption (A2) also holds, since for every $v_0 \in C^1([0,1])$, $\bar{x}_0 \in \Re$ with $cv_0'(0) = \zeta\bar{x}_0$, $v_0(0) = 0$, there exists $\bar{x} \in C^0([-r,0];\Re)$ with $\bar{x}(0) = \bar{x}_0$ that satisfies (3.10). More specifically, $\bar{x} \in C^0([-r,0];\Re)$ is given by

$$\bar{x}(-rz) = (r\zeta)^{-1}\exp(r\zeta z)v_0'(z), \text{ for } z \in [0,1] \tag{3.29}$$

for arbitrary $\eta \in C^0([-r,0];\Re^{N_1})$ with $\eta(0) = v_0(1)$. Selecting the constant function $\eta(s) = v_0(1)$ for $s \in [-r,0]$, we guarantee that inequality $\|\bar{x}\| + \|\eta\| \leq G(\|v_0\|_\infty + \|v_0'\|_\infty + |\bar{x}_0|)$ holds with $G = 1 + (r\zeta)^{-1}\exp(r\zeta)$. Using the formulas of the previous section, we can relate system (1.4), (3.27), (3.28) with output $y(t) = v(t,1)$ to the delay system

$$\begin{aligned}
\dot{x}_1(t) &= \zeta x_2(t) - \zeta x_2(t-r)\exp(-\zeta r) - \zeta x_1(t) \\
\dot{x}_2(t) &= \theta(x_2(t)) - (\mu+1)x_2(t) + \mu\zeta\int_0^1\int_{t-rl}^t x_2(s)\exp(-\zeta(t-s))dsdl + u(t) \\
y(t) &= x_1(t)
\end{aligned} \tag{3.30}$$

with initial condition that satisfies $x_1(0) = \zeta\int_{-r}^0 \exp(\zeta s)x_2(s)ds$. As mentioned above, the reactor model (1.4), (3.27), (3.28) is *not* equivalent to system (3.30): the solutions of the reactor model (1.4), (3.27), (3.28) correspond to the solutions of (3.30) only when $x_1(0) = \zeta\int_{-r}^0 \exp(\zeta s)x_2(s)ds$.

In this case, we get $x_1(t) = \zeta\int_{t-r}^t \exp(\zeta(s-t))x_2(s)ds$ for all $t \geq 0$ and consequently, the reactor model (1.4), (3.27), (3.28) turns out to be equivalent to the following system with distributed state and output delays:

$$\frac{d\bar{x}}{dt}(t) = \theta(\bar{x}(t)) - (\mu+1)\bar{x}(t) + \mu\zeta \int_0^1 \int_{t-rl}^t \bar{x}(s)\exp(-\zeta(t-s))dsdl + u(t)$$

$$y(t) = \zeta \int_{t-r}^t \exp(\zeta(s-t))\bar{x}(s)ds$$

Making use of Theorem 3.2, we can prove the following result for the reactor model.

**Theorem 3.3 (Sample-Data Observer for the Chemical Reactor):** *There exist constants $k_1, k_2, \delta, P, M, \omega > 0$ such that for every sampling sequence $\{\tau_i\}_{i=0}^\infty$ with $\tau_0 = 0$, $\lim(\tau_i) = +\infty$, $0 < \tau_{i+1} - \tau_i \leq \delta$ for $i = 0,1,...$, for every $v_0 \in C^1([0,1])$, $\bar{x}_0 \in \Re$ with $cv_0'(0) = \zeta\bar{x}_0$, $v_0(0) = 0$, $z_0 \in C^0([-r,0]; \Re^2)$ and $u \in C^0(\Re_+; \Re)$, $\xi \in L^\infty_{loc}(\Re_+; \Re)$, the solution $(\bar{x}(t), v[t], z(t), w(t))$ of (1.4), (3.27), (3.28) with*

$$\dot{z}_1(t) = \zeta z_2(t) - \zeta z_2(t-r)\exp(-\zeta r) - \zeta z_1(t) - k_1(z_1(t) - w(t)), \quad t \geq 0$$

$$\dot{z}_2(t) = \theta(z_2(t)) - (\mu+1)z_2(t) + \mu\zeta\int_0^1\int_{t-rl}^t z_2(s)\exp(-\zeta(t-s))dsdl + u(t) - k_2(z_1(t)-w(t)), \quad t \geq 0$$

$$\dot{w}(t) = \zeta z_2(t) - \zeta z_2(t-r)\exp(-\zeta r) - \zeta z_1(t), \quad t \in [\tau_i, \tau_{i+1})$$

$$\hat{x}_t = z_t, \quad t \geq 0$$

(3.31)

$$w_j(\tau_i) = \xi(\tau_i) + v(\tau_i, 1) \tag{3.32}$$

$$\hat{v}(t,z) = \zeta \int_{t-rz}^t \exp(-\zeta(t-p))\hat{x}_2(p)dp, \text{ for } t \geq 0, z \in [0,1] \tag{3.33}$$

*initial condition $(\bar{x}_0, v_0, z_0)$, corresponding to inputs $u, \xi$, exists for all $t \geq 0$ and satisfies estimate*

$$\begin{aligned}&|\bar{x}(t) - \hat{x}_2(t)| + \|v[t] - \hat{v}[t]\|_\infty \\ &\leq \exp(-\omega t)M\left(\|v_0\|_\infty + \|v_0'\|_\infty + |\bar{x}_0| + \|z_0\|\right) + P\sup_{0\leq s \leq t}\left(|\xi(s)|\exp(-\omega(t-s))\right)\end{aligned} \tag{3.34}$$

The proof of Theorem 3.3 is constructive and is given in Section 5. It should be also noticed that there is *no restriction* in the speed $c > 0$ of the transport process (or equivalently, in the delay $r > 0$).

*3.C. Stabilization of Hyperbolic PDE-ODE Loops*

For stabilization purposes, we need to assume (A1), (A2), (A3) as well as the following (stabilizability) assumption.

**(A4)** *There exists a functional $\bar{k}: \Re^{\bar{n}} \times C^0([0,1]) \to \Re^m$ such that Assumption (H4) holds for system (1.1) with $D = \{0\}$, where $f: C^0([-r,0]; \Re^n) \times \Re^m \times D \to \Re^n$ is defined by (3.12), (3.13), (3.14), (3.15) for the functional $\tilde{k}: C^0([-r,0]; \Re^n) \to \Re^m$ that satisfies the equation*

$$\tilde{k}(\bar{x},\eta,\lambda) = \bar{k}(\bar{x}_0, v_0) \tag{3.35}$$

*for every* $v_0 \in C^1([0,1])$, $\bar{x}_0 \in \Re^{\bar{n}}$ *with* $cv_0'(0) = g(0,\bar{x}_0) + \sum_{i=1}^{N_1} b_i(0)v_0(z_i) + \sum_{i=1}^{N_2} \beta_i(0)\int_0^1 \gamma_i(s)v_0(s)ds$, $v_0(0) = 0$, *where* $\bar{x} \in C^0([-r,0];\Re^{\bar{n}})$, $\eta \in C^0([-r,0];\Re^{N_1})$, $\lambda \in C^0([-r,0];\Re^{N_2})$ *with* $\bar{x}(0) = \bar{x}_0$, $\eta_i(0) = v_0(z_i)$ $(i=1,...,N_1)$, $\lambda_i(0) = \int_0^1 \gamma_i(z)v_0(z)dz$ $(i=1,...,N_2)$ *are the functions involved in (3.10).*

Under Assumption (A4), we obtain the following stabilization result.

**Corollary 3.4 (Global Stabilization of Hyperbolic PDE-ODE Loops with Observer-Based Sampled-Data Feedback):** *Suppose that Assumptions (A1), (A2), (A3), (A4) hold for system (1.2), (1.3), (1.4). Consider system (1.1) with* $D = \{0\}$, *where* $f:C^0([-r,0];\Re^n) \times \Re^m \times D \to \Re^n$, $h:C^0([-r,0];\Re^n) \to \Re^k$ *are defined by (3.12), (3.13), (3.14), (3.15), (3.16). Suppose that system (2.3) is a REO for system (1.1). Moreover, suppose that for every bounded* $S \subset C^0([-r,0];\Re^n) \times \Re^m$, *there exists a constant* $L_S \geq 0$ *such that (2.6) holds with* $R: C^0([-r,0];\Re^n) \times \Re^m \times D \to \Re^k$ *being defined by (3.20). Then there exist constants* $\delta, P, \omega > 0$ *and a function* $\bar{a} \in K_\infty$ *such that for every sampling sequence* $\{\tau_i\}_{i=0}^\infty$ *with* $\tau_0 = 0$, $\lim(\tau_i) = +\infty$, $0 < \tau_{i+1} - \tau_i \leq \delta$ *for* $i=0,1,...$, *for every* $v_0 \in C^1([0,1])$, $\bar{x}_0 \in \Re^{\bar{n}}$ *with* $cv_0'(0) = g(0,\bar{x}_0) + \sum_{i=1}^{N_1} b_i(0)v_0(z_i) + \sum_{i=1}^{N_2} \beta_i(0)\int_0^1 \gamma_i(s)v_0(s)ds$, $v_0(0) = 0$, $z_0 \in C^0([-r,0];\Re^l)$ *and* $\xi \in L^\infty_{loc}(\Re_+;\Re^k)$, *the solution* $(\bar{x}(t), v[t], z(t), w(t))$ *of (1.2), (1.3), (1.4) with (2.8), (3.21), (3.22) and*

$$u(t) = \bar{k}(H\hat{x}(\tau_i), \hat{v}[\tau_i]), \text{ for } t \in [\tau_i, \tau_{i+1}) \tag{3.36}$$

*initial condition* $(\bar{x}_0, v_0, z_0)$, *corresponding to input* $\xi$, *exists for all* $t \geq 0$ *and satisfies the following estimate:*

$$\begin{aligned} &|\bar{x}(t)| + |H\hat{x}(t)| + \|v[t]\|_\infty + \|\hat{v}[t]\|_\infty \\ &\leq \exp(-\omega t)\bar{a}\left(\|v_0\|_\infty + \|v_0'\|_\infty + |\bar{x}_0| + \|z_0\|\right) + P\sup_{0 \leq s \leq t}\left(|\xi(s)|\exp(-\omega(t-s))\right) \end{aligned} \tag{3.37}$$

*where* $H \in \Re^{\bar{n} \times n}$ *is the matrix involved in (3.24).*

The following example illustrates the use of Corollary 3.4.

**Example 3.5:** For the chemical reactor (1.4), (3.27), (3.28), Assumption (A4) holds with

$$\bar{k}(\bar{x},v) := -Q\bar{x} - \mu \int_0^1 v(z)dz, \text{ for all } x \in \Re, v \in C^0([0,1]) \tag{3.38}$$

where $Q > 0$ is a sufficiently large constant (so that $\mu + 1 + Q$ is greater than the Lipschitz constant of $\theta$). It follows from the proof of Theorem 3.3 and from Corollary 3.4 that there exist constants $\delta, P, \omega > 0$ and a function $\bar{a} \in K_\infty$ such that for every sampling sequence $\{\tau_i\}_{i=0}^\infty$ with $\tau_0 = 0$, $\lim(\tau_i) = +\infty$, $0 < \tau_{i+1} - \tau_i \leq \delta$ for $i = 0, 1, ...$, for every $v_0 \in C^1([0,1])$, $\bar{x}_0 \in \Re$ with $cv_0'(0) = \zeta \bar{x}_0$, $v_0(0) = 0$, $z_0 \in C^0([-r, 0]; \Re^2)$ and $\xi \in L_{loc}^\infty(\Re_+; \Re)$, the solution $(\bar{x}(t), v[t], z(t), w(t))$ of (1.4), (3.27), (3.28) together with (3.32), (3.31), (3.33) and

$$u(t) = -Q\hat{x}_2(\tau_i) - \mu \int_0^1 \hat{v}(\tau_i, z) dz, \text{ for } t \in [\tau_i, \tau_{i+1}) \tag{3.39}$$

with initial condition $(\bar{x}_0, v_0, z_0)$, corresponding to input $\xi$, exists for all $t \geq 0$ and satisfies estimate (3.37). ◁

## 4. Triangular Globally Lipschitz Delay Systems

In this section we consider systems of the form (1.6), where we assume that there exists a constant $L \geq 0$ such that the following inequalities hold for $i = 1, ..., n$:

$$|f_i(x_1, ..., x_i, u) - f_i(z_1, ..., z_i, u)| \leq \tilde{L} \sum_{j=1}^i \|x_j - z_j\|,$$

for all $(x_1, ..., x_i), (z_1, ..., z_i) \in C^0([-r, 0]; \Re^i)$, $u \in \Re^m$ \hfill (4.1)

Notice that systems of the form (1.6) satisfying (4.1) are Forward Complete and so satisfy Assumptions (H1), (H2) with $U = \Re^m$, $D = \Re^n$. Indeed, using (1.6), (4.1) and the triangle inequality, we get that for all $x_0 \in C^0([-r, 0]; \Re^n)$ and $u \in L_{loc}^\infty(\Re_+; \Re^m)$, $d \in L_{loc}^\infty(\Re_+; \Re^n)$ the unique solution $x$ of (1.6) with initial condition $x_0$, corresponding to inputs $u, d$, satisfies:

$$|x(t)| \leq |x_0(0)| + t \sup_{0 \leq s \leq t}(|d(s)|) + t \sup_{0 \leq s \leq t}\left(\sum_{i=1}^n |f_i(0, u(s))|\right) + n\tilde{L} \int_0^t \|x_s\| ds, \text{ for all } t \in [0, t_{max})$$

where $t_{max} \in (0, +\infty]$ is the maximal existence time of the solution. The above inequality implies the following integral inequality for all $t_1 \in (0, t_{max})$ and $t \in [0, t_1]$

$$\|x_t\| \leq \|x_0\| + t_1 \sup_{0 \leq s \leq t_1}(|d(s)|) + t_1 \sup_{0 \leq s \leq t_1}\left(\sum_{i=1}^n |f_i(0, u(s))|\right) + n\tilde{L} \int_0^t \|x_s\| ds$$

which in conjunction with the Gronwall-Bellman lemma gives:

$$\|x_t\| \leq \exp(n\tilde{L}t)\left(\|x_0\| + t_1 \sup_{0 \leq s \leq t_1}(|d(s)|) + t_1 \sup_{0 \leq s \leq t_1}\left(\sum_{i=1}^n |f_i(0, u(s))|\right)\right)$$

The above inequality and the fact that $t_1 \in (0, t_{max})$ is arbitrary guarantee (in conjunction with standard results from the theory of delay equations; see [24]) that for every $x_0 \in C^0([-r, 0]; \Re^n)$ and $u \in L^\infty_{loc}(\Re_+;\Re^m)$, $d \in L^\infty_{loc}(\Re_+;\Re^n)$, the unique solution $x(t)$ of (1.6) with initial condition $x_0$, corresponding to inputs $u, d$, exists for all $t \geq 0$ and satisfies the estimate

$$\|x_t\| \leq \exp(n\tilde{L}t)\left(\|x_0\| + t \sup_{0 \leq s \leq t}(|d(s)|) + t \sup_{0 \leq s \leq t}\left(\sum_{i=1}^n |f_i(0, u(s))|\right)\right) \quad (4.2)$$

It is also straightforward to check that Assumption (H3) holds as well with $k=1$, $L = \tilde{L}+1$, $\kappa(s) = s$ for $s \geq 0$ and $R(x, u, d) := f_1(x_1, u) + d_1$ for all $x \in C^0([-r, 0]; \Re^n)$, $(u, d) \in \Re^m \times \Re^n$.

Define the matrix $A = \{a_{i,j}: i=1,...,n, j=1,...,n\} \in \Re^{n \times n}$ by the relations

$$a_{i,i+1} = 1 \text{ for } i=1,...,n-1 \text{ and } a_{i,j} = 0 \text{ if otherwise} \quad (4.3)$$

and the vector

$$c := (1, 0, ..., 0)^T \in \Re^n \quad (4.4)$$

Since the pair of matrices $(A, c)$ is observable, there exists $K = (K_1, ..., K_n)^T \in \Re^n$ so that the matrix $(A + Kc^T) \in \Re^{n \times n}$ is Hurwitz.

Using Theorem 2.2, we are in a position to prove the following result.

**Theorem 4.1:** *There exist constants $\delta, \omega > 0$, $\theta, Q_1, Q_2, Q_3 \geq 1$ such that for every sampling sequence $\{\tau_i\}_{i=0}^\infty$ with $\tau_0 = 0$, $\lim(\tau_i) = +\infty$, $0 < \tau_{i+1} - \tau_i \leq \delta$ for all $i = 0, 1, ...$, for every $(x_0, z_0) \in C^0([-r, 0]; \Re^n) \times C^0([-r, 0]; \Re^n)$ and $u \in L^\infty_{loc}(\Re_+; \Re^m)$, $d \in L^\infty_{loc}(\Re_+; \Re^n)$, $\xi \in L^\infty_{loc}(\Re_+; \Re)$, the solution $(x(t), z(t), w(t))$ of (1.6) with*

$$\begin{aligned}
\dot{z}_i(t) &= f_i(z_{1,t}, ..., z_{i,t}, u(t)) + z_{i+1}(t) + \theta^i K_i(c^T z(t) - w(t)), i=1,...,n-1, t \geq 0 \\
\dot{z}_n(t) &= f_n(z_{1,t}, ..., z_{n,t}, u(t)) + \theta^n K_n(c^T z(t) - w(t)), t \geq 0 \\
\hat{x}_t &= z_t, t \geq 0 \\
\dot{w}(t) &= f_1(z_{1,t}, u(t)) + z_2(t), t \in [\tau_i, \tau_{i+1})
\end{aligned} \quad (4.5)$$

$$w(\tau_i) = x_1(\tau_i) + \xi(\tau_i) \quad (4.6)$$

*initial condition $(x_0, z_0)$, corresponding to inputs $u, d, \xi$, exists for all $t \geq 0$ and satisfies the following estimate*

$$\begin{aligned}
\|\hat{x}_t - x_t\| &\leq \exp(-\omega t) Q_1(\|x_0\| + \|z_0\|) \\
&+ Q_2 \sup_{0 \leq s \leq t}(|\xi(s)| \exp(-\omega(t-s))) + Q_3 \sup_{0 \leq s \leq t}(|d(s)|)
\end{aligned} \quad (4.7)$$

The proof of Theorem 4.1 is based on a combined Lyapunov analysis together with small-gain arguments. The observer (4.5), (4.6) is constructed by the combination of a high-gain observer with an inter-sample predictor. However, as the proof of Theorem 4.1 shows, the observer parameter $\theta \geq 1$ for system (1.6) has to be greater than that of the corresponding delay-free, triangular, globally Lipschitz system. More specifically, the parameter $\theta \geq 1$ depends on the maximum delay $r > 0$ (see inequality (5.43) in the proof of Theorem 4.1).

## 5. Proofs

We start with the proof of Theorem 2.2.

**Proof of Theorem 2.2:** Let $i \geq 0$ be an integer for which $z_{\tau_i}$ exists. We will show first that $z_t$ exists for all $t \in [\tau_i, \tau_{i+1}]$. Due to the regularity assumptions of Definition 2.1 and due to (2.6) there exists $t_{\max} > 0$ such that the solution of

$$\dot{z}(t) = F(z_t, w(t), u(t))$$
$$\dot{w}(t) = R(\Psi(z_t), u(t), 0) \tag{5.1}$$

is defined on $t \in [\tau_i, t_{\max})$, where $t_{\max} > \tau_i$ is the maximal existence time of the solution of (5.1). If $t_{\max} > \tau_{i+1}$ then there is nothing to show. Consequently, we next focus on the case $t_{\max} \leq \tau_{i+1}$. Define:

$$v(t) := w(t) - h(x_t) \tag{5.2}$$

Notice that the component of the solution $z_t$ of (5.1) is actually a solution of

$$\dot{z}(t) = F(z_t, h(x_t) + v(t), u(t)) \tag{5.3}$$

for all $t \in [0, t_{\max})$. Therefore, by virtue of (2.5) and since $\omega \leq \sigma$, the following estimate holds for all $t \in [0, t_{\max})$:

$$\|\hat{x}_t - x_t\| \leq \exp(-\omega t) a\left(\|x_0\| + \|z_0\|\right) + \gamma \sup_{0 \leq s \leq t} \left(|v(s)| \exp(-\omega(t-s))\right) + \sup_{0 \leq s \leq t} \left(g(|d(s)|)\right) \tag{5.4}$$

Using (2.1), (2.2), (2.8) and (2.9), we get for all $t \in [0, t_{\max})$:

$$|v(t)| \leq |\xi(q(t))| + L \int_{q(t)}^{t} \|\hat{x}_s - x_s\| ds + (t - q(t)) \sup_{q(t) \leq s \leq t} \left(\kappa(|d(s)|)\right) \tag{5.5}$$

where

$$q(t) = \max\{\tau_i : \tau_i \leq t\}. \tag{5.6}$$

Using (5.5) and the fact that $t - q(t) \leq \delta$ (a consequence of definition (5.6) and the fact that $0 < \tau_{i+1} - \tau_i \leq \delta$ for $i = 0, 1, \ldots$), we obtain for all $t \in [0, t_{\max})$:

$$|v(t)| \exp(\omega t) \leq |\xi(q(t))| \exp(\omega q(t)) \exp(\omega(t - q(t)))$$
$$+ L \exp(\omega t) \int_{q(t)}^{t} \exp(-\omega s) \|\hat{x}_s - x_s\| \exp(\omega s) ds + \delta \exp(\omega t) \sup_{0 \leq s \leq t}\left(\kappa(|d(s)|)\right)$$
$$\leq |\xi(q(t))| \exp(\omega q(t)) \exp(\omega \delta)$$
$$+ L \frac{\exp(\omega \delta) - 1}{\omega} \sup_{q(t) \leq s \leq t} \left(\|\hat{x}_s - x_s\| \exp(\omega s)\right) + \delta \exp(\omega t) \sup_{0 \leq s \leq t}\left(\kappa(|d(s)|)\right)$$

which directly implies the following estimate for all $t \in [0, t_{\max})$:

$$|v(t)|\exp(\omega t) \leq \exp(\sigma \delta) \sup_{0 \leq s \leq t}\left(|\xi(s)|\exp(\omega s)\right)$$
$$+ L \frac{\exp(\omega \delta) - 1}{\omega} \sup_{0 \leq s \leq t}\left(\|\hat{x}_s - x_s\|\exp(\omega s)\right) + \delta \exp(\omega t) \sup_{0 \leq s \leq t}\left(\kappa(|d(s)|)\right) \quad (5.7)$$

Combining estimates (5.4) and (5.7) we get for all $t \in [0, t_{\max})$:

$$\|\hat{x}_t - x_t\|\exp(\omega t) \leq a\left(\|x_0\| + \|z_0\|\right) + \gamma \exp(\omega \delta) \sup_{0 \leq s \leq t}\left(|\xi(s)|\exp(\omega s)\right)$$
$$+ \gamma L \frac{\exp(\omega \delta) - 1}{\omega} \sup_{0 \leq s \leq t}\left(\|\hat{x}_s - x_s\|\exp(\omega s)\right) \quad (5.8)$$
$$+ \exp(\omega t)\left(\gamma \delta \sup_{0 \leq s \leq t}\left(\kappa(|d(s)|)\right) + \sup_{0 \leq s \leq t}\left(g(|d(s)|)\right)\right)$$

Since $B := \gamma L \frac{\exp(\sigma \delta) - 1}{\sigma} < 1$, it follows from (5.8) and definition $\tilde{g}(s) := (1 - B)^{-1}(g(s) + \gamma \delta \kappa(s))$ that estimate (2.10) holds for all $t \in [0, t_{\max})$.

Estimates (5.7) and (2.10) show that $|v(t)|$ is bounded on $[0, t_{\max})$. Extending $v(t)$ on the interval $[t_{\max}, +\infty)$ in a way that $|v(t)|$ is bounded on $\Re_+$, it follows that the solution of (5.3) exists for $t \in [\tau_i, t_{\max}]$, which shows that $t_{\max} > \tau_i$ is not the maximal existence time of the solution of (5.1), a contradiction.

Therefore, the case $t_{\max} \leq \tau_{i+1}$ cannot arise and consequently, $t_{\max} > \tau_i$. By induction and since $\lim(\tau_i) = +\infty$, it follows that $z_t$ exists for all $t \geq 0$. Moreover, estimate (2.10) holds for all $t \geq 0$. The proof is complete. ◁

We next continue with the proof of Corollary 2.4.

**Proof of Corollary 2.4:** First of all, Assumptions (H1), (H2), (H3), (H4) guarantee that the closed-loop system (1.1) with (2.8), (2.9) and (2.13) is forward complete, i.e., its unique solution exists for all $t \geq 0$, for arbitrary initial conditions, arbitrary inputs and arbitrary sampling sequences $\{\tau_i\}_{i=0}^{\infty}$ with $\tau_0 = 0$, $\lim(\tau_i) = +\infty$, $0 < \tau_{i+1} - \tau_i \leq \delta$ for $i = 0, 1, \ldots$, where $\delta > 0$ is sufficiently small (satisfies the conditions of Theorem 2.2).

Using a slight modification in the proof of Theorem 5 in [36] and applying a small-gain analysis instead of Halanay's inequality, we are in a position to show that Assumption (H4) guarantees the existence of constants $\tilde{\delta}, \bar{\gamma}, \bar{M}, \omega > 0$ and a function $\bar{g} \in K$ such that for every sampling sequence $\{\tau_i\}_{i=0}^{\infty}$ with $\tau_0 = 0$, $\lim(\tau_i) = +\infty$, $0 < \tau_{i+1} - \tau_i \leq \tilde{\delta}$ for $i = 0, 1, \ldots$, for every $x_0 \in C^0([-r, 0]; \Re^n)$ and $\tilde{d} \in L^{\infty}_{loc}(\Re_+; \Re^m)$, $d \in L^{\infty}_{loc}(\Re_+; D)$, the solution $x(t)$ of

$$\dot{x}(t) = f(x_t, u(t), d(t))$$
$$u(t) = \tilde{k}\left(x_{\tau_i}\right) + \tilde{d}(\tau_i), t \in [\tau_i, \tau_{i+1}) \quad (5.9)$$

with initial condition $x_0$, corresponding to inputs $\tilde{d}, d$, exists for all $t \geq 0$ and satisfies the following estimate

$$\|x_t\| \leq \exp(-\omega t)\bar{M}\|x_0\| + \bar{\gamma}\sup_{0\leq s\leq t}\left(|\tilde{d}(s)|\exp(-\omega(t-s))\right) + \sup_{0\leq s\leq t}\left(\bar{g}(|d(s)|)\right), \forall t \geq 0 \quad (5.10)$$

Notice that the component $x(t)$ of the solution of the closed-loop system (1.1) with (2.8), (2.9) and (2.13) coincides with the solution of (5.9) with the same initial condition, same input $d$ and input $\tilde{d} \in L^\infty_{loc}(\mathfrak{R}_+; \mathfrak{R}^m)$ defined by

$$\tilde{d}(t) = \tilde{k}(\hat{x}_{\tau_i}) - \tilde{k}(x_{\tau_i}), \text{ for } t \in [\tau_i, \tau_{i+1}) \text{ and } i \geq 0 \quad (5.11)$$

It follows from Assumption (H4) and (5.11) that the following inequality holds:

$$|\tilde{d}(t)| \leq \bar{L}\|\hat{x}_{\tau_i} - x_{\tau_i}\|, \text{ for } t \in [\tau_i, \tau_{i+1}) \text{ and } i \geq 0 \quad (5.12)$$

Selecting $\bar{\delta} > 0$ sufficiently small so that both estimates (2.10) and (5.10) hold for sampling sequences $\{\tau_i\}_{i=0}^\infty$ with $\tau_0 = 0$, $\lim(\tau_i) = +\infty$, $0 < \tau_{i+1} - \tau_i \leq \bar{\delta}$, for $i = 0,1,...$, we get from (2.10) and (5.12) that the following estimate holds for all $t \geq 0$:

$$\begin{aligned}|\tilde{d}(t)|\exp(\omega t) &\leq \bar{L}(1-B)^{-1}\exp(\omega\bar{\delta})a\left(\|x_0\| + \|z_0\|\right) \\ &+ \bar{L}(1-B)^{-1}\gamma\exp(2\omega\bar{\delta})\sup_{0\leq s\leq t}\left(|\xi(s)|\exp(\omega s)\right) + \bar{L}\exp(\omega t)\sup_{0\leq s\leq t}\left(\tilde{g}(|d(s)|)\right)\end{aligned} \quad (5.13)$$

where $\tilde{g}(s) := (1-B)^{-1}(g(s) + \gamma\bar{\delta}\kappa(s))$ and $B := \gamma L \dfrac{\exp(\omega\bar{\delta}) - 1}{\omega} < 1$. Combining (2.10), (5.10), (5.13) and using the triangle inequality $\|\hat{x}_t\| \leq \|x_t\| + \|\hat{x}_t - x_t\|$, we obtain (2.14) for appropriate constant $\hat{\gamma} > 0$ and appropriate functions $\hat{g}, \hat{a} \in K$. Notice that if $a \in K_\infty$ (the function involved in (2.5), (2.10) and (5.13)) is linear then $\hat{a}$ is linear. The proof is complete. ◁

The proof of Theorem 3.2 follows.

**Proof of Theorem 3.2:** Under Assumptions (A1), (A2), (A3), the solution $v \in C^1(\mathfrak{R}_+ \times [0,1])$ and $\bar{x} \in C^1(\mathfrak{R}_+; \mathfrak{R}^{\bar{n}})$ of system (1.2), (1.3), (1.4), (1.5) with initial conditions $v_0 \in C^1([0,1])$, $\bar{x}_0 \in \mathfrak{R}^{\bar{n}}$ with $cv'_0(0) = g(0,\bar{x}_0) + \sum_{i=1}^{N_1} b_i(0)v_0(z_i) + \sum_{i=1}^{N_2} \beta_i(0)\int_0^1 \gamma_i(s)v_0(s)ds$, $v_0(0) = 0$, corresponding to input $u \in C^0(\mathfrak{R}_+; \mathfrak{R}^m)$, is expressed by (3.4), where $x(t) = (\bar{x}(t), \eta(t), \lambda(t))$ is the solution of system (1.1) with $D = \{0\}$, $f : C^0([-r,0]; \mathfrak{R}^n) \times \mathfrak{R}^m \times D \to \mathfrak{R}^n$, $h : C^0([-r,0]; \mathfrak{R}^n) \to \mathfrak{R}^k$ being defined by (3.12), (3.13), (3.14), (3.15), (3.16), initial condition provided by Assumption (A2) and corresponding to the same input $u \in C^0(\mathfrak{R}_+; \mathfrak{R}^m)$.

It follows that Theorem 3.2 is a direct consequence of Theorem 2.2, Lemma 3.1 and the fact that there exists a constant $G > 0$ such that the initial condition of system (1.1) satisfies the inequality $\|x\| \leq G(\|v_0\|_\infty + \|v'_0\|_\infty + |\bar{x}_0|)$ (recall Assumption (A2) and definition (3.11)). Indeed, using (3.4), (3.6), (3.11), (3.22) and the fact that there exists a constant $L>0$ such that the inequality $\max_{0\leq z\leq 1}\left(|g(z,x)-g(z,y)|\right)\leq L|x-y|$ holds (recall Assumption (A1)), we obtain the existence of a constant $B>0$ (that depends only on $c>0$, $a\in C^0([0,1])$, $b_i\in C^1([0,1])$ ($i=1,...,N_1$), $\gamma_i\in C^0([0,1]\times[0,1])$, $\beta_i\in C^1([0,1])$ ($i=1,...,N_2$) and $g\in C^1([0,1]\times\Re^n;\Re)$) such that the following inequality holds for all $t\geq 0$:

$$\|\hat{v}[t]-v[t]\|_\infty \leq B\|\hat{x}_t - x_t\| \tag{5.14}$$

Inequality (3.23) is obtained by combining (2.10), (5.14) and using the fact the initial condition of system (1.1) satisfies the inequality $\|x\|\leq G\left(\|v_0\|_\infty + \|v_0'\|_\infty + |\bar{x}_0|\right)$. The proof is complete. ◁

We next provide the proof of Theorem 3.3.

**Proof of Theorem 3.3:** By virtue of Theorem 3.2, it suffices to prove that for appropriate selection of the constants $k_1, k_2 \in \Re$, the system

$$\dot{z}_1(t) = \zeta z_2(t) - \zeta z_2(t-r)\exp(-\zeta r) - \zeta z_1(t) - k_1(z_1(t)-y(t))$$

$$\dot{z}_2(t) = \theta(z_2(t)) - (\mu+1)z_2(t) + \mu\zeta\int_0^1\int_{t-rl}^t z_2(s)\exp(-\zeta(t-s))dsdl + u(t) - k_2(z_1(t)-y(t)) \tag{5.15}$$

$$\hat{x}_t = z_{2,t}$$

is a REO for system (3.30). In order to prove this fact, we consider the Lyapunov functional

$$V(t) = \frac{R}{2}(z_1(t)-x_1(t))^2 + Q\int_{t-r}^t (z_2(s)-x_2(s))^2 \exp(-\zeta(t-s))ds$$

$$+ \frac{1}{2}(z_2(t)-x_2(t)-bz_1(t)+bx_1(t))^2 \tag{5.16}$$

where $R, b, Q > 0$ are constants to be selected. For every $\xi \in L^\infty_{loc}(\Re_+;\Re)$, the time derivative of $V(t)$ along the trajectories of (5.15), (3.30) with $y(t) = x_1(t) + \xi(t)$ is given by

$$\begin{aligned}\dot{V}(t) &= -\left((k_1+\zeta)R - R\zeta b - Qb^2\right)E_1^2(t) \\ &\quad -(\mu+1+b\zeta-Q)E_2^2(t) - \zeta Q\int_{t-r}^t (z_2(s)-x_2(s))^2\exp(-\zeta(t-s))ds \\ &\quad +\left(R\zeta + 2bQ + b(k_1+\zeta) - k_2 - (\mu+1+b\zeta)b\right)E_1(t)E_2(t) \\ &\quad +\left(bE_2(t) - RE_1(t)\right)\zeta\exp(-\zeta r)(z_2(t-r)-x_2(t-r)) \\ &\quad -Q(z_2(t-r)-x_2(t-r))^2\exp(-\zeta r) \\ &\quad +(\theta(z_2(t))-\theta(x_2(t)))E_2(t) + \left((k_2-bk_1)E_2(t) + Rk_1E_1(t)\right)\xi(t) \\ &\quad +\mu\zeta E_2(t)\int_0^1\int_{t-rl}^t \exp(-\zeta(t-s))(z_2(s)-x_2(s))dsdl\end{aligned} \tag{5.17}$$

for $t\geq 0$ a.e., where

$$E_1(t) := z_1(t) - x_1(t) \tag{5.18}$$

$$E_2(t) := z_2(t) - x_2(t) - bE_1(t) \tag{5.19}$$

Since $\theta : \Re \to \Re$ is a globally Lipschitz function, there exists a constant $\Phi > 0$ such that

$$|\theta(z_2(t)) - \theta(x_2(t))| \leq \Phi |z_2(t) - x_2(t)| \leq \Phi |E_2(t)| + \Phi b |E_1(t)| \tag{5.20}$$

In the above inequality we have used the triangle inequality and the definitions (5.18), (5.19). Consequently, by using the (Young) inequalities

$$\Phi b |E_1(t)||E_2(t)| \leq \Phi E_2^2(t) + \Phi \frac{b^2}{4} E_1^2(t)$$

$$E_2(t)(z_2(t-r) - x(t-r)) \leq \frac{1}{2} E_2^2(t) + \frac{1}{2}(z_2(t-r) - x(t-r))^2$$

$$|E_1(t)||z_2(t-r) - x(t-r)| \leq \frac{1}{2} E_1^2(t) + \frac{1}{2}(z_2(t-r) - x(t-r))^2$$

$$|k_2 - bk_1||E_2(t)||\xi(t)| \leq E_2^2(t) + \frac{|k_2 - bk_1|^2}{4} \xi^2(t)$$

$$|k_1||E_1(t)||\xi(t)| \leq E_1^2(t) + \frac{k_1^2}{4} \xi^2(t)$$

we obtain from (5.17), (5.20) the following inequality for $t \geq 0$ a.e.:

$$\begin{aligned}
\dot{V}(t) \leq &-\left( (k_1 + \zeta)R - R\zeta b - Qb^2 - \Phi \frac{b^2}{4} - \frac{1}{2} R\zeta \exp(-\zeta r) - R \right) E_1^2(t) \\
&-\left( \mu + b\zeta - Q - 2\Phi - \frac{1}{2} b\zeta \exp(-\zeta r) \right) E_2^2(t) \\
&+ \left( R\zeta + 2bQ + b(k_1 + \zeta) - k_2 - (\mu + 1 + b\zeta)b \right) E_1(t) E_2(t) \\
&+ \left( \frac{1}{2}(b+R)\zeta - Q \right) \exp(-\zeta r)(z_2(t-r) - x(t-r))^2 \\
&- \zeta Q \int_{t-r}^{t} (z_2(s) - x(s))^2 \exp(-\zeta(t-s)) ds + \frac{1}{4}\left( |k_2 - bk_1|^2 + Rk_1^2 \right) \xi^2(t) \\
&+ \mu \zeta E_2(t) \int_0^1 \int_{t-rl}^{t} \exp(-\zeta(t-s))(z_2(s) - x_2(s)) ds \, dl
\end{aligned} \tag{5.21}$$

Using the Cauchy-Schwarz inequality (twice), we bound the double integral appearing in the right hand side of (5.21) as follows:

$$\left| \int_0^1 \int_{t-rl}^t \exp(-\zeta(t-s))(z_2(s)-x_2(s))\,ds\,dl \right| \leq \int_0^1 \left| \int_{t-rl}^t \exp(-\zeta(t-s))(z_2(s)-x_2(s))\,ds \right| dl$$

$$\leq \int_0^1 \int_{t-rl}^t \exp(-\zeta(t-s))|z_2(s)-x_2(s)|\,ds\,dl$$

$$\leq \left( \int_0^1 \left( \int_{t-rl}^t \exp(-\zeta(t-s))|z_2(s)-x_2(s)|\,ds \right)^2 dl \right)^{1/2}$$

$$\leq \left( \int_0^1 \left( \int_{t-r}^t \exp(-\zeta(t-s))|z_2(s)-x_2(s)|\,ds \right)^2 dl \right)^{1/2} \leq \int_{t-r}^t \exp(-\zeta(t-s))|z_2(s)-x_2(s)|\,ds$$

$$\leq \left( \int_{t-r}^t \exp(-\zeta(t-s))\,ds \right)^{1/2} \left( \int_{t-r}^t \exp(-\zeta(t-s))|z_2(s)-x_2(s)|^2\,ds \right)^{1/2}$$

$$\leq \left( \frac{1-\exp(-\zeta r)}{\zeta} \right)^{1/2} \left( \int_{t-r}^t \exp(-\zeta(t-s))|z_2(s)-x_2(s)|^2\,ds \right)^{1/2}$$

Using the above inequality in conjunction with estimate (5.21), we obtain the following inequality for $t \geq 0$ a.e.:

$$\dot{V}(t) \leq -\left( (k_1+\zeta)R - R\zeta b - Qb^2 - \Phi\frac{b^2}{4} - \frac{1}{2}R\zeta\exp(-\zeta r) - R \right)E_1^2(t)$$

$$-\left( \mu + b\zeta - Q - 2\Phi - \frac{1}{2}b\zeta\exp(-\zeta r) \right)E_2^2(t)$$

$$+\left( R\zeta + 2bQ + b(k_1+\zeta) - k_2 - (\mu+1+b\zeta)b \right)E_1(t)E_2(t)$$

$$+\left( \frac{1}{2}(b+R)\zeta - Q \right)\exp(-\zeta r)(z_2(t-r)-x_2(t-r))^2 \quad (5.22)$$

$$-\zeta Q \int_{t-r}^t (z_2(s)-x_2(s))^2 \exp(-\zeta(t-s))\,ds + \frac{1}{4}\left( |k_2-bk_1|^2 + Rk_1^2 \right)\xi^2(t)$$

$$+\mu\zeta|E_2(t)|\left( \frac{1-\exp(-\zeta r)}{\zeta} \right)^{1/2} \left( \int_{t-r}^t (z_2(s)-x_2(s))^2 \exp(-\zeta(t-s))\,ds \right)^{1/2}$$

Finally, using the Young inequality

$$|E_2(t)|\left( \frac{1-\exp(-\zeta r)}{\zeta} \right)^{1/2} \left( \int_{t-r}^t (z_2(s)-x_2(s))^2 \exp(-\zeta(t-s))\,ds \right)^{1/2}$$

$$\leq \frac{1}{2}E_2^2(t) + \frac{1-\exp(-\zeta r)}{2\zeta} \int_{t-r}^t (z_2(s)-x_2(s))^2 \exp(-\zeta(t-s))\,ds$$

in conjunction with estimate (5.22), we obtain the following inequality for $t \geq 0$ a.e.:

$$\dot{V}(t) \leq -\left((k_1 + \zeta)R - R\zeta b - Qb^2 - \Phi\frac{b^2}{4} - \frac{1}{2}R\zeta\exp(-\zeta r) - R\right)E_1^2(t)$$

$$-\left(\mu + b\zeta - Q - 2\Phi - \frac{1}{2}b\zeta\exp(-\zeta r) - \frac{1}{2}\mu\zeta\right)E_2^2(t)$$

$$+\left(R\zeta + 2bQ + b(k_1 + \zeta) - k_2 - (\mu + 1 + b\zeta)b\right)E_1(t)E_2(t) \tag{5.23}$$

$$-\left(Q - \frac{1}{2}(b+R)\zeta\right)\exp(-\zeta r)(z_2(t-r) - x_2(t-r))^2 + \frac{1}{4}\left(|k_2 - bk_1|^2 + Rk_1^2\right)\xi^2(t)$$

$$-\left(\zeta Q - \mu\frac{1-\exp(-\zeta r)}{2}\right)\int_{t-r}^{t}(z_2(s) - x_2(s))^2\exp(-\zeta(t-s))ds$$

By selecting

$$R = \frac{2\mu(1-\exp(-\zeta r))}{\zeta^2}, \quad b = \frac{4\Phi + (\mu + R + 1)\zeta}{\zeta(1-\exp(-\zeta r))} \tag{5.24}$$

$$Q = \frac{1}{2}(b+R)\zeta, \quad k_2 = R\zeta + b(b+R)\zeta + b(k_1+\zeta) - (\mu+1+b\zeta)b \tag{5.25}$$

$$k_1 = \zeta b + \frac{1}{2R}(b+R)b^2\zeta + \Phi\frac{b^2}{4R} + \frac{1}{2}\zeta\exp(-\zeta r) + 1 \tag{5.26}$$

we obtain from (5.23) and (5.16) the following differential inequality for $t \geq 0$ a.e.:

$$\dot{V}(t) \leq -\frac{\zeta}{2}V(t) + \frac{1}{4}\left(|k_2 - bk_1|^2 + Rk_1^2\right)\xi^2(t) \tag{5.27}$$

Applying Lemma 2.12 in [24] in conjunction with (5.27), we get for all $t \geq 0$:

$$V(t) \leq \exp\left(-\frac{\zeta}{2}t\right)V(0) + \frac{1}{4}\left(|k_2 - bk_1|^2 + Rk_1^2\right)\int_0^t \exp\left(-\frac{\zeta}{2}(t-s)\right)\xi^2(s)ds \tag{5.28}$$

Notice that the quadratic form $S(x) = \frac{R}{2}x_1^2 + \frac{1}{2}(x_2 - bx_1)^2$ on $\Re^2$ is positive definite. Consequently, there exists a constant $K_1 > 0$ such that $K_1|x|^2 \leq S(x)$ for all $x \in \Re^2$. Using this fact in conjunction with definition (5.16) and bounding the integral in the right hand side of (5.28) in the following way for any $\sigma \in \left(0, \frac{\zeta}{4}\right)$

$$\int_0^t \exp\left(-\frac{\zeta}{2}(t-s)\right)|\xi(s)|^2 ds$$

$$\leq \int_0^t \exp\left(-\frac{\zeta}{2}(t-s)\right)\exp(-2\sigma s)ds \sup_{0\leq s\leq t}\left(|\xi(s)|^2\exp(2\sigma s)\right)$$

$$\leq 2\frac{\exp(-2\sigma t) - \exp\left(-\frac{\zeta}{2}t\right)}{\zeta - 4\sigma}\sup_{0\leq s\leq t}\left(|\xi(s)|^2\exp(2\sigma s)\right)$$

$$\leq \frac{2\exp(-2\sigma t)}{\zeta - 4\sigma}\sup_{0\leq s\leq t}\left(|\xi(s)|^2\exp(2\sigma s)\right)$$

we obtain from (5.28) for all $t \geq 0$:

$$|z(t)-x(t)|^2 \leq K_1^{-1}\exp(-2\sigma t)V(0) + \frac{|k_2-bk_1|^2 + Rk_1^2}{2K_1(\zeta-4\sigma)}\sup_{0\leq s\leq t}\left(\xi^2(s)\exp(-2\sigma(t-s))\right) \quad (5.29)$$

Definition (5.16) implies that there exists a constant $K_2 \geq K_1$ (independent of $z_t, x_t$) such that $K_2\left(\|x_t\|^2 + \|z_t\|^2\right) \geq V(t)$. Therefore, we obtain from (5.29) for all $t \geq 0$:

$$|z(t)-x(t)| \leq \sqrt{\frac{K_2}{K_1}}\exp(-\sigma t)(\|x_0\|+\|z_0\|) + \sqrt{\frac{|k_2-bk_1|^2 + Rk_1^2}{2K_1(\zeta-4\sigma)}}\sup_{0\leq s\leq t}\left(|\xi(s)|\exp(-\sigma(t-s))\right) \quad (5.30)$$

Notice that due to the fact that $K_2 \geq K_1$, inequality (5.30) actually holds for all $t \geq -r$. Consequently, we obtain from (5.30) for all $t \geq -r$:

$$\sup_{-r\leq s\leq t}\left(|z(s)-x(s)|\exp(\sigma s)\right) \leq \sqrt{\frac{K_2}{K_1}}(\|x_0\|+\|z_0\|) + \sqrt{\frac{|k_2-bk_1|^2 + Rk_1^2}{2K_1(\zeta-4\sigma)}}\sup_{0\leq s\leq t}\left(|\xi(s)|\exp(\sigma s)\right) \quad (5.31)$$

Using the fact that $\sup_{t-r\leq s\leq t}\left(|z(s)-x(s)|\exp(\sigma s)\right) \geq \exp(\sigma(t-r))\|z_t-x_t\|$, we obtain from (5.31) for all $t \geq 0$:

$$\|z_t-x_t\| \leq \sqrt{\frac{K_2}{K_1}}\exp(-\sigma(t-r))(\|x_0\|+\|z_0\|) + \sqrt{\frac{|k_2-bk_1|^2 + Rk_1^2}{2K_1(\zeta-4\sigma)}}\exp(\sigma r)\sup_{0\leq s\leq t}\left(|\xi(s)|\exp(-\sigma(t-s))\right)$$

(5.32)

Estimate (5.32) holds for all $\xi \in L^\infty_{loc}(\Re_+;\Re)$ and for all $\sigma \in \left(0, \frac{\zeta}{4}\right)$ and shows that (5.15) is a REO for system (3.30). The proof is complete. ◁

The proof of Corollary 3.4 follows.

**Proof of Corollary 3.4:** Under Assumptions (A1), (A2), (A3), (A4), the solution $(x(t), v[t], z(t), w(t))$ of (1.2), (1.3), (1.4) with (2.8), (3.21), (3.22) and (3.36) is expressed by means (3.4), (3.7), (3.8) using the solution $(x(t), z(t), w(t))$ of system (1.1) with (2.8), (2.9) and (2.13) with $D = \{0\}$, where $f: C^0([-r,0];\Re^n) \times \Re^m \times D \to \Re^n$, $h: C^0([-r,0];\Re^n) \to \Re^k$ are defined by (3.12), (3.13), (3.14), (3.15), (3.16), with initial condition provided by Assumption (A2).

Following the same procedure in the proof of Theorem 3.2, we conclude that Corollary 3.4 is a direct consequence of Corollary 2.4 and the fact that there exists a constant $G > 0$ such that the initial condition of system (1.1) satisfies the inequality $\|x\| \leq G(\|v_0\|_\infty + \|v'_0\|_\infty + |x_0|)$ (recall Assumption (A2)). The proof is complete. ◁

Finally, we provide the proof of Theorem 4.1.

**Proof of Theorem 4.1:** By virtue of Theorem 2.2, it suffices to show that the following system

$$\dot{z}_i(t) = f_i(z_{1,t},...,z_{i,t},u(t)) + z_{i+1}(t) + \theta^i K_i(c^T z(t) - y(t)), \; i = 1,...,n-1$$
$$\dot{z}_n(t) = f_n(z_{1,t},...,z_{n,t},u(t)) + \theta^n K_n(c^T z(t) - y(t)) \quad (5.33)$$
$$\hat{x}_t = z_t$$

is a ROE for system (1.6).

The proof is based on the quadratic error Lyapunov function $V(t) := e^T(t)\Delta_\theta^{-1} P \Delta_\theta^{-1} e(t)$, where $e(t) := z(t) - x(t)$, $\Delta_\theta := diag(\theta, \theta^2,...,\theta^n)$ and $P \in \Re^{n \times n}$ is a symmetric positive definite matrix that satisfies $P(A + Kc^T) + (A^T + cK^T)P + 2\mu I \leq 0$ for certain constant $\mu > 0$. Notice that the following inequalities hold:

$$\theta^{-i}\left|f_i(x_{1,t},...,x_{i,t},u(t)) - f_i(z_{1,t},...,z_{i,t},u(t))\right| \leq n\tilde{L}\|\varepsilon_t\|, \text{ for } i=1,...,n \text{ and } t \geq 0, \quad (5.34)$$

where

$$\varepsilon(t) = \Delta_\theta^{-1} e(t) \quad (5.35)$$

Indeed, inequalities (5.34) follow from (4.1), the fact that $\theta \geq 1$, the fact that $\sum_{j=1}^{i}\|\varepsilon_{j,t}\| \leq n\|\varepsilon_t\|$ and definition (5.35). For every $u \in L_{loc}^\infty(\Re_+;\Re^m)$, $d \in L_{loc}^\infty(\Re_+;\Re^n)$, $\xi \in L_{loc}^\infty(\Re_+;\Re)$ we obtain from (1.6) and (5.33) with $y(t) = x_1(t) + \xi(t)$:

$$\dot{e}(t) = (A + \Delta_\theta Kc^T)e(t) + g(x_t,e_t,u(t)) - d(t) - \Delta_\theta K\xi(t), \text{ for } t \geq 0 \text{ a.e.} \quad (5.36)$$

where $g(x_t,e_t,u(t)) = (f_1(x_{1,t}+e_{1,t},u(t)) - f_1(x_{1,t},u(t)),...,f_n(x_t+e_t,u(t)) - f_n(x_t,u(t)))^T$. Therefore, using the identities $\Delta_\theta^{-1}A = \theta A \Delta_\theta^{-1}$, $c^T = \theta c^T \Delta_\theta^{-1}$ and definition (5.35), we get for $t \geq 0$ a.e.:

$$\dot{V}(t) \leq -2\theta\mu|\varepsilon(t)|^2 + 2|\varepsilon(t)||P||\Delta_\theta^{-1}g(x_t,e_t,u(t))| + 2|\varepsilon(t)||PK||\xi(t)| + 2|\varepsilon(t)||P\Delta_\theta^{-1}||d(t)| \quad (5.37)$$

Using (5.34) and (5.37), the fact that $|\Delta_\theta^{-1}| \leq \theta^{-1}$ (a consequence of definition $\Delta_\theta := diag(\theta,\theta^2,...,\theta^n)$ and the fact that $\theta \geq 1$) in conjunction with the inequalities

$$2|\varepsilon(t)||PK||\xi(t)| \leq \frac{\theta\mu}{3}|\varepsilon(t)|^2 + \frac{3}{\theta\mu}|PK|^2|\xi(t)|^2, \quad 2|\varepsilon(t)||P\Delta_\theta^{-1}||d(t)| \leq \frac{\theta\mu}{3}|\varepsilon(t)|^2 + \frac{3}{\theta\mu}|P\Delta_\theta^{-1}|^2|d(t)|^2,$$

$$2n\sqrt{n}|\varepsilon(t)||P|\tilde{L}\|\varepsilon_t\| \leq \frac{\theta\mu}{3}|\varepsilon(t)|^2 + \frac{3n^3}{\theta\mu}|P|^2\tilde{L}^2\|\varepsilon_t\|^2, \text{ we obtain for } t \geq 0 \text{ a.e.:}$$

$$\dot{V}(t) \leq -\theta\mu|\varepsilon(t)|^2 + \frac{3n^3}{\theta\mu}|P|^2\tilde{L}^2\|\varepsilon_t\|^2 + \frac{3}{\theta\mu}|PK|^2|\xi(t)|^2 + \frac{3}{\theta^3\mu}|P|^2|d(t)|^2 \quad (5.38)$$

Setting $\varphi := \frac{\theta\mu}{2|P|}$ and using the fact that $V(t) = \varepsilon^T(t)P\varepsilon(t) \leq |P||\varepsilon(t)|^2$, we obtain from (5.38) for $t \geq 0$ a.e.:

$$\dot{V}(t) \leq -2\varphi V(t) + \frac{3n^3}{\theta\mu}|P|^2\tilde{L}^2\|\varepsilon_t\|^2 + \frac{3}{\theta\mu}|PK|^2|\xi(t)|^2 + \frac{3}{\theta^3\mu}|P|^2|d(t)|^2 \quad (5.39)$$

Applying Lemma 2.12 in [24] in conjunction with (5.39), we get for all $t \geq 0$:

$$V(t) \leq \exp(-2\varphi t)V(0) + \frac{3n^3}{\theta\mu}|P|^2 \tilde{L}^2 \int_0^t \exp(-2\varphi(t-s))\|\varepsilon_s\|^2 \, ds$$
$$+ \frac{3}{\theta\mu}|PK|^2 \int_0^t \exp(-2\varphi(t-s))|\xi(s)|^2 \, ds + \frac{3}{\theta^3\mu}|P|^2 \int_0^t \exp(-2\varphi(t-s))|d(s)|^2 \, ds \quad (5.40)$$

Since $P \in \mathfrak{R}^{n \times n}$ is positive definite there exists a constant $c_1 > 0$ such that $c_1|x|^2 \leq x^T P x$ for all $x \in \mathfrak{R}^n$. Using this fact, the fact that $V(t) = \varepsilon^T(t) P \varepsilon(t) \leq |P||\varepsilon(t)|^2$ and bounding the three integrals in the right hand side of (5.40) in the following way (similarly for the other two integrals) for any $\sigma \in (0, \varphi)$

$$\int_0^t \exp(-2\varphi(t-s))|\xi(s)|^2 \, ds$$
$$\leq \int_0^t \exp(-2\varphi(t-s))\exp(-2\sigma s) \, ds \sup_{0 \leq s \leq t}\left(|\xi(s)|^2 \exp(2\sigma s)\right)$$
$$\leq \frac{\exp(-2\sigma t) - \exp(-2\varphi t)}{2(\varphi - \sigma)} \sup_{0 \leq s \leq t}\left(|\xi(s)|^2 \exp(2\sigma s)\right)$$
$$\leq \frac{\exp(-2\sigma t)}{2(\varphi - \sigma)} \sup_{0 \leq s \leq t}\left(|\xi(s)|^2 \exp(2\sigma s)\right)$$

we obtain from (5.40) for all $t \geq 0$:

$$|\varepsilon(t)|^2 \exp(2\sigma t) \leq \frac{|P|}{c_1}|\varepsilon(0)|^2 + \frac{3n^3|P|^2 \tilde{L}^2}{2\theta\mu c_1(\varphi - \sigma)} \sup_{0 \leq s \leq t}\left(\|\varepsilon_s\|^2 \exp(2\sigma s)\right)$$
$$+ \frac{3|PK|^2}{2\theta\mu c_1(\varphi - \sigma)} \sup_{0 \leq s \leq t}\left(|\xi(s)|^2 \exp(2\sigma s)\right) + \frac{3|P|^2}{2\theta^3\mu c_1(\varphi - \sigma)} \sup_{0 \leq s \leq t}\left(|d(s)|^2 \exp(2\sigma s)\right) \quad (5.41)$$

Using the fact that $\sup_{t-r \leq s \leq t}\left(|\varepsilon(s)|^2 \exp(2\sigma s)\right) \geq \exp(2\sigma(t-r)) \sup_{t-r \leq s \leq t}\left(|\varepsilon(s)|^2\right)$, we obtain from (5.41) the following estimate for all $t \geq 0$:

$$|\varepsilon(t)|\exp(\sigma t) \leq \sqrt{\frac{|P|}{c_1}}|\varepsilon(0)| + \sqrt{\frac{3n^3|P|^2 \tilde{L}^2 \exp(2\sigma r)}{2\theta\mu c_1(\varphi - \sigma)}} \sup_{-r \leq s \leq t}\left(|\varepsilon(s)|\exp(\sigma s)\right)$$
$$+ \sqrt{\frac{3|PK|^2}{2\theta\mu c_1(\varphi - \sigma)}} \sup_{0 \leq s \leq t}\left(|\xi(s)|\exp(\sigma s)\right) + \sqrt{\frac{3|P|^2}{2\theta^3\mu c_1(\varphi - \sigma)}} \sup_{0 \leq s \leq t}\left(|d(s)|\exp(\sigma s)\right) \quad (5.42)$$

Selecting $\theta \geq 1$ so large so that

$$\frac{3n^3|P|^2 \tilde{L}^2 \exp(\varphi r)}{\theta\mu c_1 \varphi} < 1 \quad (5.43)$$

and using the fact that $\sup_{-r \leq s \leq t}\left(|\varepsilon(s)|\exp(\sigma s)\right) \leq \sup_{-r \leq s \leq 0}\left(|\varepsilon(s)|\exp(\sigma s)\right) + \sup_{0 \leq s \leq t}\left(|\varepsilon(s)|\exp(\sigma s)\right)$ we obtain from (5.42) for $\sigma = \varphi/2$ and for all $t \geq 0$:

$$\sup_{0\leq s\leq t}\left(|\varepsilon(s)|\exp(\sigma s)\right)\leq (1-\Omega)^{-1}\sqrt{\frac{|P|}{c_1}}|\varepsilon(0)|+(1-\Omega)^{-1}\Omega\sup_{-r\leq s\leq 0}\left(|\varepsilon(s)|\exp(\sigma s)\right)$$
$$+(1-\Omega)^{-1}\left(\sqrt{\frac{3|PK|^2}{\theta\mu c_1\varphi}}\sup_{0\leq s\leq t}\left(|\xi(s)|\exp(\sigma s)\right)+\sqrt{\frac{3|P|^2}{\theta^3\mu c_1\varphi}}\sup_{0\leq s\leq t}\left(|d(s)|\exp(\sigma s)\right)\right) \quad (5.44)$$

where $\Omega:=\sqrt{\dfrac{3n^3|P|^2\tilde{L}^2\exp(\varphi r)}{\theta\mu c_1\varphi}}<1$. Due to the fact that $|P|\geq c_1$, it follows from (5.44) that the following estimate holds for all $t\geq -r$:

$$\sup_{-r\leq s\leq t}\left(|\varepsilon(s)|\exp(\sigma s)\right)\leq (1-\Omega)^{-1}\left(\sqrt{\frac{|P|}{c_1}}+\Omega\right)\|\varepsilon_0\|$$
$$+(1-\Omega)^{-1}\left(\sqrt{\frac{3|PK|^2}{\theta\mu c_1\varphi}}\sup_{0\leq s\leq t}\left(|\xi(s)|\exp(\sigma s)\right)+\sqrt{\frac{3|P|^2}{\theta^3\mu c_1\varphi}}\sup_{0\leq s\leq t}\left(|d(s)|\exp(\sigma s)\right)\right) \quad (5.45)$$

Using the fact that $\sup_{t-r\leq s\leq t}\left(|\varepsilon(s)|\exp(\sigma s)\right)\geq \exp(\sigma(t-r))\|\varepsilon_t\|$, we obtain from (5.45) the following estimate for all $t\geq 0$:

$$\|\varepsilon_t\|\leq (1-\Omega)^{-1}\left(\sqrt{\frac{|P|}{c_1}}+\Omega\right)\|\varepsilon_0\|\exp(-\sigma(t-r))$$
$$+(1-\Omega)^{-1}\sqrt{\frac{3|PK|^2}{\theta\mu c_1\varphi}}\exp(\sigma r)\sup_{0\leq s\leq t}\left(|\xi(s)|\exp(-\sigma(t-s))\right) \quad (5.46)$$
$$+(1-\Omega)^{-1}\sqrt{\frac{3|P|^2}{\theta^3\mu c_1\varphi}}\exp(\sigma r)\sup_{0\leq s\leq t}\left(|d(s)|\right)$$

Finally, definition (5.35) and the fact that $\theta\geq 1$ imply the inequalities $|e(t)|\leq \theta^n|\varepsilon(t)|$ and $|\varepsilon(t)|\leq \theta^{-1}|e(t)|$. Consequently, we get from (5.46):

$$\|e_t\|\leq (1-\Omega)^{-1}\left(\sqrt{\frac{|P|}{c_1}}+\Omega\right)\theta^{n-1}\|e_0\|\exp(-\sigma(t-r))$$
$$+(1-\Omega)^{-1}\theta^n\sqrt{\frac{3|PK|^2}{\theta\mu c_1\varphi}}\exp(\sigma r)\sup_{0\leq s\leq t}\left(|\xi(s)|\exp(-\sigma(t-s))\right) \quad (5.47)$$
$$+(1-\Omega)^{-1}\theta^{n-1}\sqrt{\frac{3|P|^2}{\theta\mu c_1\varphi}}\exp(\sigma r)\sup_{0\leq s\leq t}\left(|d(s)|\right)$$

Since $e_t=\hat{x}_t-x_t$ and $\|e_0\|\leq \|x_0\|+\|z_0\|$, inequality (5.47) shows that (5.33) is a ROE for system (1.6). The proof is complete. ◁

# 6. Concluding Remarks

The present work studied the problem of designing sampled-data observers and observer-based, sampled-data, output feedback stabilizers for systems with both discrete and distributed, state and output time-delays. The obtained results were applied to time delay systems of strict-feedback structure, transport PDEs with nonlocal terms, and feedback interconnections of ODEs with a transport PDE. The study constitutes a unified theoretical framework to deal with (sampled-data observer design and output-feedback stabilization for) ODEs, delay systems and PDEs, while these were generally dealt with separately.

There are many things to be done in the case of ODE-PDE loops. It is important to extend the results to the case of more than one transport PDE. This extension is a topic for future research. The results of the present work can be easily extended to the case of networked control systems with uniformly globally exponentially protocols (see for instance [6]).